\documentclass[11pt,a4paper]{article}
\usepackage[T1]{fontenc}
\usepackage[utf8]{inputenc}
\usepackage{lmodern,amsmath,amssymb,amsthm,mathtools}
\usepackage[margin=28mm]{geometry}
\usepackage{microtype,enumitem,booktabs,array}
\usepackage{xcolor}
\usepackage[colorlinks=true,linkcolor=blue!45!black,citecolor=blue!45!black,urlcolor=blue!45!black]{hyperref}
\hypersetup{pdftitle={A discrete Baras--Goldstein theorem and the sharp lattice Hardy threshold},pdfauthor={Carlos Lizama}}
\numberwithin{equation}{section}
\newtheorem{theorem}{Theorem}[section]
\newtheorem{proposition}[theorem]{Proposition}
\newtheorem{lemma}[theorem]{Lemma}
\newtheorem{corollary}[theorem]{Corollary}
\theoremstyle{definition}
\theoremstyle{remark}\newtheorem{remark}[theorem]{Remark}
\newcommand{\Z}{\mathbb Z}
\newcommand{\R}{\mathbb R}
\newcommand{\E}{\mathcal E_N}
\newcommand{\M}{M_N}
\newcommand{\K}{\mathcal K_N}
\newcommand{\T}{T_N}
\newcommand{\Om}{\Omega_N}
\newcommand{\ip}[2]{\langle #1,#2\rangle}
\newcommand{\norm}[1]{\lVert #1\rVert}
\newcommand{\supp}{\operatorname{supp}}
\newcommand{\Dom}{\operatorname{Dom}}
\newcommand{\Adj}{\operatorname{Adj}}

\newcommand{\dd}{\,\mathrm d}

\setlist[enumerate]{label=\textup{(\roman*)},leftmargin=*}
\allowdisplaybreaks[2]
\title{A discrete Baras--Goldstein theorem\\and the sharp lattice Hardy threshold}
\author{Carlos Lizama\\[4pt]
\small Departamento de Matem\'atica y Ciencia de la Computaci\'on\\
\small Universidad de Santiago de Chile, Santiago, Chile\\[2pt]
\small \href{mailto:carlos.lizama@usach.cl}{carlos.lizama@usach.cl}\\
\small ORCID: \href{https://orcid.org/0000-0002-9807-1100}{0000-0002-9807-1100}}
\date{}
\begin{document}
\maketitle
\begin{abstract}
We establish a lattice counterpart of the Baras--Goldstein phenomenon for the heat equation with an inverse-square potential. On punctured integer lattices, the sharp discrete Hardy constant determines a threshold independent of the mesh size. At or below this threshold the heat semigroups are uniformly contractive in the natural discrete $L^2$ spaces. Above it, solutions with fixed nonnegative, nonzero macroscopic initial data diverge as the mesh tends to zero. We determine their logarithmic growth rate uniformly on bounded observation regions and for times bounded away from zero. Below the discrete threshold, the cell extensions converge in $L^2$ to the continuum Friedrichs evolution, uniformly for all times. In dimensions $N\ge9$, a nonempty interval of coupling constants supports a contractive continuum evolution while the sampled lattice solutions diverge under mesh refinement. We express the discrete threshold as the reciprocal norm of an explicit weighted Green operator, with a kernel given by modified Bessel functions and the Dirichlet correction at the origin. Finally, an explicit rational function approximates the threshold with absolute error less than $10^9N^{-5}$ for every integer dimension $N\ge50$. These bounds give explicit sufficient conditions for stability and blow-up.
\end{abstract}


\section{Introduction and main results}

The inverse-square potential has the same scaling as the Laplacian. This balance makes the sharp Hardy constant a boundary between two different behaviors of heat evolution. For the continuum equation
\[
 \partial_tu=\Delta u+\frac{\lambda}{|x|^2}u,
 \qquad N\ge3,
\]
the critical value is $A_N=(N-2)^2/4$. Baras and Goldstein \cite{BG84}, answering questions posed by Brezis and Lions, proved existence for suitable nonnegative data when $0\le\lambda\le A_N$ and nonexistence of nontrivial nonnegative solutions when $\lambda>A_N$. Their formulation through bounded approximations also identifies complete instantaneous blow-up: in the supercritical regime the approximating solutions tend to infinity at every positive time and every observation point. The conditions on the initial data are part of the existence theorem; see \cite[Proposition~2.1 and Theorem~2.2]{BG84}.

This connection between a functional inequality and positive evolution persists in broader settings. Goldstein and Zhang treated the Heisenberg group \cite{GZ01} and variable elliptic coefficients \cite{GZ03}. Higher-order operators were studied by Galaktionov and Kamotski \cite{GK10}. Fractional Laplacians with Hardy potentials were investigated by Bogdan, Grzywny, Jakubowski and Pilarczyk \cite{BGJP19}; singular drift terms were considered by Okazawa, Sobajima and Yokota \cite{OSY14} and, in the fractional setting, by Kinzebulatov, Sem\"enov and Szczypkowski \cite{KSS21}. The recent proof of Goldstein, Goldstein, K\"ombe and Rhandi \cite{GGKR24} emphasizes the role of scaling in producing instantaneous blow-up.

The present paper studies the corresponding problem on $\Z^N\setminus\{0\}$ with nearest-neighbor interactions. For each fixed mesh, the discrete Laplacian and the inverse-square potential are bounded operators. Their heat evolution is therefore globally defined. The critical behavior is expressed by a sharp stability threshold on each lattice and by blow-up under mesh refinement. We show how the sharp discrete Hardy constant controls both the stability of the lattice evolution and its approximation of the continuum problem. Our results connect this spectral threshold to a precise macroscopic growth law, establish strong convergence to the continuum Friedrichs evolution below the threshold, and make the threshold accessible through an explicit Green representation and a quantitative rational approximation.

Hardy inequalities themselves have a long discrete history; see Kufner, Maligranda and Persson \cite{KMP06}. The modern theory on graphs constructs optimal spatial weights from positive superharmonic functions and Green functions \cite{KPP18,Fischer24}. Keller and Lemm \cite{KL23} analyze the inverse-square behavior of such weights on Euclidean lattices, including operators with variable coefficients. Fischer and Rose \cite{FischerRose26} develop related Poincar\'e--Hardy inequalities on manifolds, homogeneous trees and weakly spherically symmetric graphs. We consider the prescribed Euclidean weight $|x|^{-2}$, characterize its sharp scalar coefficient and obtain a quantitative approximation in high dimension. The optimality of a spatially varying Hardy weight and the sharpness of a scalar coefficient for a fixed weight are distinct notions.

\subsection{The lattice equation and its threshold}

Fix an integer $N\ge3$ and let $\Om=\Z^N\setminus\{0\}$. We use the Euclidean norm $|n|$ and denote the coordinate vectors by $e_1,\ldots,e_N$. Every function on $\Om$ is extended by zero at the origin. The space $c_{00}(\Om)$ consists of the complex-valued functions on $\Om$ with finite support. We take inner products to be linear in the first argument. On $\ell^2(\Om)$ with counting measure, the positive Dirichlet Laplacian is
\begin{equation}\label{eq:laplacian}
 L_Nv(n)=2Nv(n)-\sum_{j=1}^N\bigl(v(n+e_j)+v(n-e_j)\bigr),
 \qquad n\in\Om.
\end{equation}
Its diagonal remains $2N$ at the neighbors of the origin, so edges incident to the origin retain their Dirichlet energy. For $v\in c_{00}(\Om)$, define the energy and the weighted mass by
\begin{equation}\label{eq:energy-mass}
 \E[v]=\sum_{n\in\Z^N}\sum_{j=1}^N|v(n+e_j)-v(n)|^2,
 \qquad
 \M[v]=\sum_{n\in\Om}\frac{|v(n)|^2}{|n|^2}.
\end{equation}
Each unoriented edge is counted once. The sharp Hardy constant for this energy and weight is
\begin{equation}\label{eq:constant}
 C(N)=\inf_{0\ne v\in c_{00}(\Om)}\frac{\E[v]}{\M[v]}.
\end{equation}
Proposition~\ref{prop:coercivity} proves that $0<C(N)\le2N$.

The dimensional size of this constant has been studied through several complementary methods. Gupta \cite[Corollary~1.5]{Gupta23} proved $C(N)\asymp N$, where $\asymp$ denotes two-sided comparability. Gupta \cite[Theorem~1]{Gupta26} and, independently, Huang and Ye \cite[Theorem~1.1]{HuangYe26} then identified the leading coefficient: $C(N)/(2N)\to1$. Their results also cover higher-order Hardy--Rellich inequalities.

The exact low-dimensional values were obtained by Alpay for dimension three \cite{Alpay3} and dimensions five through eight \cite[Theorem~2.1]{AlpayLow}, and by Lizama for dimensions three and four \cite{LizamaLowDimensions}. The latter work also gives strict upper bounds from dimension nine. Together these results yield the transition, already stated in \cite[Corollary~2.2]{AlpayLow},
\begin{equation}\label{eq:known-transition}
 C(N)=A_N\quad\text{for }3\le N\le8,
 \qquad C(N)<A_N\quad\text{for }N\ge9,
\end{equation}
where $A_N=(N-2)^2/4$ is the continuum Hardy constant. Thus the discrete and continuum thresholds first differ in dimension nine. The subsequent total-defect and spectral analysis in \cite{LizamaTransition} describes how the interaction of the vertices nearest the origin with the exterior produces this change. Here we study its consequence for heat evolution and obtain an explicit representation and quantitative approximation of the threshold.

For $h>0$, equip $h\Om$ with the measure $h^N$. The corresponding Hilbert space $\ell_h^2=\ell^2(h\Om,h^N)$ has norm
\[
 \|u\|_{\ell_h^2}^2=h^N\sum_{n\in\Om}|u(hn)|^2.
\]
For $v(n)=u(hn)$, define $\Delta_hu(hn)=-h^{-2}L_Nv(n)$ and $V_hu(x)=|x|^{-2}u(x)$. Both operators are bounded on $\ell_h^2$. For $\lambda\ge0$, the heat equation is
\begin{equation}\label{eq:heat-h}
 \partial_tu_h(t,x)=\Delta_hu_h(t,x)+\frac{\lambda}{|x|^2}u_h(t,x),
 \qquad x\in h\Om,\quad t>0.
\end{equation}
The potential is a prescribed function of $x$, so this is the linear Cauchy problem
\[
 u_h'(t)=-H^h_{\lambda,N}u_h(t),\qquad
 H^h_{\lambda,N}=-\Delta_h-\lambda V_h.
\]
Its bounded generator $-H^h_{\lambda,N}$ defines the uniformly continuous semigroup below. For each initial datum $u_{0,h}\in\ell_h^2$, the unique global solution is
\[
 S_{\lambda,h}(t)=e^{-tH^h_{\lambda,N}},\qquad
 u_h(t)=S_{\lambda,h}(t)u_{0,h}.
\]
On the unit lattice, write $V=V_1$, $H_{\lambda,N}=L_N-\lambda V$ and $\mu_N(\lambda)=-\inf\sigma(H_{\lambda,N})$.

Our first theorem gives the discrete Baras--Goldstein threshold and the growth rate above it.
\begin{theorem}\label{thm:main-bg}
Let $N\ge3$. The sharp Hardy constant on $h\Om$ equals $C(N)$ for every $h>0$, and for every $t\ge0$,
\begin{equation}\label{eq:main-growth}
 \|S_{\lambda,h}(t)\|_{\ell_h^2\to\ell_h^2}
 =\begin{cases}
 1,&0\le\lambda\le C(N),\\
 e^{\mu_N(\lambda)t/h^2},&\lambda>C(N),
 \end{cases}
\end{equation}
with $\mu_N(\lambda)>0$ in the second case. If $\lambda>C(N)$ and $f$ is a nonzero, nonnegative, continuous function with compact support in $\R^N\setminus\{0\}$, the solution with $u_h(0,hn)=f(hn)$ satisfies
\begin{equation}\label{eq:main-log-limit}
 \lim_{h\downarrow0}\ 
 \sup_{\substack{t\ge\tau\\n\in\Om,\ |hn|\le R}}
 \bigl|h^2\log u_h(t,hn)-\mu_N(\lambda)t\bigr|=0
 \qquad(R,\tau>0).
\end{equation}
\end{theorem}

In particular, the supercritical solutions tend to $+\infty$ uniformly for $t\ge\tau$ and $0<|hn|\le R$. Proposition~\ref{bg:pointwise-rate} gives a quantitative estimate. At $\lambda=A_N$, the norm identity recovers \cite[Corollary~3.17]{LizamaTransition}; the present formulation extends that statement to arbitrary couplings, with threshold $C(N)$. The new parabolic statement identifies how the negative spectral bottom governs divergence for fixed macroscopic data throughout bounded observation regions, including points whose lattice coordinates escape to infinity as the mesh is refined.

The threshold is independent of $h$ because the energy and the weighted mass both acquire the factor $h^{N-2}$. The proof of pointwise divergence follows the diagonal-growth strategy of \cite[Theorem~3.1]{GGKR24}, followed by propagation along lattice paths. General logarithmic heat-kernel asymptotics for fixed vertices are available in \cite[Theorems~2.2 and~3.1]{KLVW15}. Here time, sampled data and observation vertices vary simultaneously with $h$; the path estimate gives the required uniform control.

\subsection{Relation to the continuum evolution}

The stable range raises a different question: whether the discrete solutions approximate the continuum solution. The functional framework, well-posedness and large-time behavior of the continuum equation were developed by V\'azquez and Zuazua \cite{VZ00}. For $0\le\lambda\le A_N$, let $\mathsf H_{\lambda,N}^{F}$ be the Friedrichs realization of $-\Delta-\lambda|x|^{-2}$ on $L^2(\R^N)$. Extend lattice functions constantly over cells by $I_h$, with value zero on the cell at the origin, and let $P_h=I_h^*$ take cell averages. Section~\ref{cont:section} proves that, for $0\le\lambda<C(N)$ and every $f\in L^2(\R^N)$,
\begin{equation}\label{eq:intro-continuum}
 \sup_{t\ge0}\bigl\|I_hS_{\lambda,h}(t)P_hf
              -e^{-t\mathsf H_{\lambda,N}^{F}}f\bigr\|_{L^2}
 \longrightarrow0\qquad(h\downarrow0).
\end{equation}
The uniformity over the entire time half-line combines approximation at finite times with control of the long-time evolution. The proof establishes convergence of the resolvents by using the strict Hardy bound to control discrete gradients, and then applies the general relation between resolvent and semigroup convergence \cite{KS03}. At $\lambda=C(N)$, the same conclusion holds if $C(N)=A_N$ or
$C(N)\le A_N-1$; see Proposition~\ref{cont:endpoint-esa} and
Theorem~\ref{thm:endpoint-critical-mosco}. These endpoint arguments use
uniqueness of the relevant nonnegative continuum realization in place of
strict coercivity. Together these results include $\lambda=C(N)$ for $3\le N\le8$ and $N\ge10$. The identification of the continuum limit at $N=9$, $\lambda=C(9)$, is not resolved here.

For $N\ge9$, the interval $C(N)<\lambda\le A_N$ has a different behavior: the continuum Friedrichs evolution is contractive, whereas the sampled lattice solutions diverge by Theorem~\ref{thm:main-bg}. The dimensional transition therefore becomes a macroscopic obstruction to convergence for this discretization. Refining the mesh preserves the discrete threshold and cannot remove the instability in this interval. Thus the singularity makes the lattice stability threshold decisive even though the finite differences are consistent away from the origin.

This distinction also separates the present approximation from conforming finite element approximations. Della Pietra et al.\ \cite{DPFIMPZ24} and Ignat and Zuazua \cite{IZ26} approximate the continuum Hardy quotient with its exact integrals; the latter obtain an excess of order $|\log h|^{-2}$ above $A_N$ on bounded smooth convex domains containing the origin, with shape-regular, quasiuniform conforming meshes. Their critical forms remain nonnegative. Our nearest-neighbor energy and sampled inverse-square mass instead have the mesh-independent sharp ratio $C(N)$. Standard lattice-to-continuum norm resolvent results such as \cite{NT21}, which assume a continuous potential bounded from below, do not cover this attractive singular potential.

\subsection{An exact Green representation}

Let $I_k$ denote the modified Bessel function of the first kind of order $k$. For $z\in\Z^N$, define
\begin{equation}\label{eq:g}
 g_N(z)=\int_0^\infty e^{-2Nt}\prod_{j=1}^N I_{|z_j|}(2t)\dd t,
\end{equation}
For $n,m\in\Om$, set
\begin{equation}\label{eq:GD}
 G_N^D(n,m)=g_N(n-m)-\frac{g_N(n)g_N(m)}{g_N(0)},
\end{equation}
and
\begin{equation}\label{eq:Tkernel}
 T_N(n,m)=\frac{G_N^D(n,m)}{|n||m|}.
\end{equation}
The Bessel heat kernel is classical \cite{Ciaurri,DLMF}. The Dirichlet correction is the single-vertex form of the killed Green identity \cite{LawlerLimic,KLW21}.

The zero-energy Birman--Schwinger principle yields the following explicit representation of the sharp threshold. This concrete weighted inverse links the Dirichlet Green kernel to bilateral variational bounds and to the supercritical growth exponent.
\begin{theorem}\label{thm:green}\label{thm:main}
For every $N\ge3$, the kernel \eqref{eq:Tkernel} defines a bounded positive selfadjoint operator $T_N$ on $\ell^2(\Om)$, with strictly positive entries, and
\begin{equation}\label{eq:mainnorm}
 C(N)=\frac1{\|T_N\|}.
\end{equation}
\end{theorem}

Proposition~\ref{gr:envelopes} turns this identity into exact upper and lower variational bounds. Section~\ref{gr:shifted} also characterizes the growth exponent $\mu_N(\lambda)$ through the norm of the weighted resolvent $V^{1/2}(L_N+s)^{-1}V^{1/2}$ for $s>0$.

\subsection{A rational approximation with an explicit error}

Set
\begin{equation}\label{eff:R-intro}
 \begin{gathered}
 R(N)=2N-4-\frac{20}{3}
 \frac{N-c}{(N-a)(N-c)-b},\\
 a=\frac{253}{60},\qquad b=\frac{10609}{1800},\qquad
 c=\frac{33448197}{5941040}.
 \end{gathered}
\end{equation}
The known expansion $C(N)=2N-4-20/(3N)+O(N^{-2})$ was proved in \cite[Theorem~1.1]{LizamaAsymptotics}. Appendix~\ref{eff:asym-appendix} recovers those terms and extends the expansion through $N^{-5}$. The parameters in $R(N)$ match the coefficients through $N^{-4}$. The next theorem makes this approximation effective for the actual infinite-lattice constant: both the error bound and its range of validity are explicit.

\begin{theorem}\label{eff:main}
For every integer $N\ge50$,
\begin{equation}\label{eff:main-bound}
 |C(N)-R(N)|<\frac{10^9}{N^5}.
\end{equation}
\end{theorem}

Consequently, a finite rational expression attains any prescribed absolute accuracy $\varepsilon>0$ in every integer dimension $N\ge N_\varepsilon$, where
\begin{equation}\label{eff:epsilon-threshold}
 N_\varepsilon=\max\left\{50,
 \left\lceil(10^9/\varepsilon)^{1/5}\right\rceil\right\}.
\end{equation}

The coefficient $10^9$ is an explicit unoptimized majorant. It gives an error less than $10^{-6}$ for $N\ge1000$ and less than $10^{-11}$ for $N\ge10000$; at $N=50$ the bound is $3.2$. The proof uses an exterior Hardy profile, a scalar Schur complement, a finite exact identity for the resolvent, and a holomorphic estimate for finitely many return moments. Appendix~\ref{eff:asym-appendix} supplies an analytic finite-core derivation of the coefficient identities used in the cancellation. The resulting error bound yields the explicit dimensional threshold in \eqref{eff:epsilon-threshold}.

The paper proceeds from the energy space to the evolution theorem and the continuum limit, then to the Green representation and its bilateral bounds. An explicit exterior estimate separates the first lattice sphere from the rest of the weighted Hardy operator and leads to the rational approximation. The final section translates the approximation into sufficient stability and blow-up conditions. The appendix contains the detailed coefficient calculation.

\section{Energy space and the weighted Hardy operator}\label{sec:forms}\label{gr:forms}

\subsection{A positive Hardy lower bound}

We first establish $C(N)>0$ to justify the inverse operator used below.
The standard ground-state identity on graphs \cite{KPP18,KLW21} converts a positive supersolution into a Hardy lower bound. The proof below fixes our edge normalization.

\begin{lemma}\label{lem:gs}
Let $h(x)>0$ on $\Om$ and $h(0)=0$. For every $u\in c_{00}(\Om)$,
\begin{equation}\label{eq:gs}
 \E[u]=
 \sum_{\substack{\{x,y\}\subset\Om\\x\sim y}}
 h(x)h(y)\left|\frac{u(x)}{h(x)}-\frac{u(y)}{h(y)}\right|^2
 +\sum_{x\in\Om}\frac{L_Nh(x)}{h(x)}|u(x)|^2.
\end{equation}
The first sum counts every unoriented edge once.
\end{lemma}
\begin{proof}
For an edge with both endpoints in $\Om$, expansion gives
\begin{align*}
 &|u(x)-u(y)|^2
 -h(x)h(y)\left|\frac{u(x)}{h(x)}-\frac{u(y)}{h(y)}\right|^2\\
 &\hspace{12mm}=
 \left(1-\frac{h(y)}{h(x)}\right)|u(x)|^2
 +\left(1-\frac{h(x)}{h(y)}\right)|u(y)|^2.
\end{align*}
An edge from $x\in\Om$ to the origin contributes $|u(x)|^2$ to the energy. Adding these contributions yields the coefficient
$2N-\sum_{y\sim x,\,y\ne0}h(y)/h(x)=L_Nh(x)/h(x)$.
Only finitely many edges meet the support of $u$, so all rearrangements are finite.
\end{proof}

Positive Hardy bounds on this lattice are already known, for example from
\cite{Gupta23}. The elementary proof below supplies the coercivity needed
for our inverse construction, without using its sharp value.

\begin{proposition}\label{prop:coercivity}
For every $N\ge3$ there is a constant $c_N>0$ such that
\begin{equation}\label{eq:coercivity}
 c_N\M[u]\le\E[u]\qquad(u\in c_{00}(\Om)).
\end{equation}
In particular, $0<C(N)\le2N$.
\end{proposition}
\begin{proof}
Write $p=(N-2)/2>0$ and $H_a(x)=(a+|x|^2)^{-p/2}$ on $\R^N$. Set $S=a+|x|^2$. Direct differentiation gives
\begin{equation}\label{eq:contshift}
 -\Delta H_a(x)=p\bigl[Na+(N-p-2)|x|^2\bigr]S^{-p/2-2}
 \ge p^2S^{-p/2-1}.
\end{equation}
For each coordinate,
\begin{align*}
 \partial_j^4H_a={}&3p(p+2)S^{-p/2-2}
 -6p(p+2)(p+4)x_j^2S^{-p/2-3}\\
 &+p(p+2)(p+4)(p+6)x_j^4S^{-p/2-4}.
\end{align*}
Consequently
\[
 |\partial_j^4H_a(x)|\le B_pS^{-p/2-2},\qquad
 B_p=p(p+2)\bigl[3+6(p+4)+(p+4)(p+6)\bigr].
\]
If $a\ge2$ and $|t|\le1$, then
$a+|x+te_j|^2\ge(a+|x|^2)/2$. Taylor's formula with remainder therefore gives
\begin{equation}\label{eq:taylorcoercive}
 \left|L_N^{\rm full}H_a(x)+\Delta H_a(x)\right|
 \le D_N S^{-p/2-2},\qquad
 D_N=\frac{NB_p}{12}\,2^{p/2+2}.
\end{equation}
Choose $a\ge\max\{2,2D_N/p^2\}$. Equations \eqref{eq:contshift} and \eqref{eq:taylorcoercive} imply
\[
 \frac{L_N^{\rm full}H_a(x)}{H_a(x)}
 \ge\frac{p^2}{S}-\frac{D_N}{S^2}
 \ge\frac{p^2}{2S}.
\]
Now let $h=H_a$ on $\Om$ and set $h(0)=0$. This modification increases the Laplacian at the neighbors of zero and leaves it unchanged elsewhere. Since $|x|\ge1$ on $\Om$,
\[
 \frac{L_Nh(x)}{h(x)}\ge\frac{p^2}{2(a+|x|^2)}
 \ge\frac{p^2}{2(a+1)}\frac1{|x|^2}.
\]
Lemma~\ref{lem:gs} proves \eqref{eq:coercivity} with $c_N=p^2/[2(a+1)]$. The trial $u=\mathbf1_{\{e_1\}}$ has energy $2N$ and mass one, giving the upper bound.
\end{proof}

\subsection{Closure of the form}

Let $D$ denote multiplication by $|x|$ on $\ell^2(\Om)$. On $c_{00}(\Om)$ put
\begin{equation}\label{eq:qform}
 \mathfrak q_N[f]=\E[Df].
\end{equation}
Let $\nabla_d$ be the difference operator indexed by the edges $(x,x+e_j)$, including edges incident with the origin. The operator $A_0f=\nabla_d(Df)$ is closable: if $f_n\to0$ in $\ell^2$ and $A_0f_n\to v$ in the edge space, then each coordinate of $v$ is zero, because each edge difference depends on two fixed coordinates of $f_n$. Thus $v=0$.

We write $\mathfrak q_N$ for the closed form $\norm{\overline A_0f}_2^2$ and $\mathcal Q_N$ for its domain. Equivalently, $\mathcal Q_N$ is the completion of $c_{00}(\Om)$ in the form norm. Proposition~\ref{prop:coercivity} gives
\begin{equation}\label{eq:qpositive}
 \mathfrak q_N[f]\ge c_N\norm f_2^2\qquad(f\in\mathcal Q_N).
\end{equation}
Let $\K$ be the positive selfadjoint operator associated with this closed form. On finitely supported functions it acts as $DL_ND$. Its domain is characterized by
\[
 \Dom(\K)=\{f\in\mathcal Q_N:\ \exists g\in\ell^2,
 \ \mathfrak q_N(f,v)=\ip gv\text{ for every }v\in\mathcal Q_N\},
 \qquad \K f=g.
\]
Since $D$ maps $c_{00}(\Om)$ bijectively onto itself, the form-core property gives
\begin{equation}\label{eq:bottom}
 \inf\sigma(\K)
 =\inf_{0\ne f\in\mathcal Q_N}\frac{\mathfrak q_N[f]}{\norm f_2^2}
 =C(N).
\end{equation}
In particular, $\K^{-1}$ is bounded. The inverse constructed below therefore acts in the weighted Hardy space; the ordinary operator $L_N$ has spectral bottom zero.

A cutoff argument places finite-energy physical profiles in the closed form
domain, where they can be used in the variational formulas below.

\begin{lemma}\label{lem:cutoff}
If $u(0)=0$, $\E[u]<\infty$ and $\M[u]<\infty$, then $f=D^{-1}u$ belongs to $\mathcal Q_N$ and $\mathfrak q_N[f]=\E[u]$.
\end{lemma}
\begin{proof}
Choose $\chi_R(x)=\eta(|x|/R)$, where $\eta$ is Lipschitz, equals one on $[0,1]$ and zero on $[2,\infty)$. Then $u_R=\chi_Ru$ has finite support and $\M[u-u_R]\to0$. On an edge $(x,y)$,
\[
 ((1-\chi_R)u)(y)-((1-\chi_R)u)(x)
 =(1-\chi_R(y))(u(y)-u(x))+u(x)(\chi_R(x)-\chi_R(y)).
\]
The square of the first term sums to a tail of the finite energy. The second contributes at most
\[
 \frac{c}{R^2}\sum_{R-1\le|x|\le2R+1}|u(x)|^2
 \le c'\sum_{|x|\ge R-1}\frac{|u(x)|^2}{|x|^2},
\]
which tends to zero. Hence $D^{-1}u_R\to D^{-1}u$ in the form norm, and the assertion follows from closure.
\end{proof}

\section{The heat equation and refinement of the mesh}\label{bg:section}

We now prove the parabolic assertion in Theorem~\ref{thm:main-bg}.
Throughout this section, $N\ge3$, $\Om=\Z^N\setminus\{0\}$, and functions
on $\Om$ have value zero at the origin. We write
\[
 (Vv)(n)=|n|^{-2}v(n),\qquad
 H_{\lambda,N}=L_N-\lambda V,\qquad \lambda\ge0.
\]
The positivity of $C(N)$ is established in the Hardy preliminaries. The
arguments below use its variational definition and do not require its
numerical value.

\subsection{Exact scaling and the Hardy threshold}

For $h>0$, set $\Om^h=h\Z^N\setminus\{0\}$ and define
\[
 \ell_h^2(\Om^h)=
 \left\{u:\Om^h\longrightarrow\mathbb C:
       h^N\sum_{x\in\Om^h}|u(x)|^2<\infty\right\},
 \qquad
 \norm{u}_{2,h}^2=h^N\sum_{x\in\Om^h}|u(x)|^2.
\]
The discrete Laplacian and the attractive Schr\"odinger operator are
\begin{align}
 (\Delta_hu)(x)
 &=\frac1{h^2}\sum_{j=1}^N
       \bigl(u(x+he_j)+u(x-he_j)-2u(x)\bigr),\label{bg:mesh-laplacian}\\
 H_{\lambda,N}^hu
 &=-\Delta_hu-\lambda|x|^{-2}u.\label{bg:mesh-H}
\end{align}
Here and below $u(0)=0$ in every term involving the origin. In particular,
the diagonal of $-\Delta_h$ is $2N/h^2$ also at its neighboring vertices.
The evolution considered in this section is
\begin{equation}\label{bg:mesh-evolution}
 \partial_tu_h+H_{\lambda,N}^hu_h=0,
 \qquad u_h(0)=u_{0,h}.
\end{equation}

We reduce all meshes to $h=1$ before studying stability. At $\lambda=A_N$, this recovers the scaling in \cite[Corollary~3.17]{LizamaTransition}; the same calculation applies to every coupling.

\begin{lemma}\label{bg:scaling}
The map
\[
 \mathcal U_h:\ell_h^2(\Om^h)\longrightarrow\ell^2(\Om),
 \qquad (\mathcal U_hu)(n)=h^{N/2}u(hn),
\]
is unitary and satisfies
\begin{equation}\label{bg:unitary-scaling}
 \mathcal U_hH_{\lambda,N}^h\mathcal U_h^{-1}
       =h^{-2}H_{\lambda,N}.
\end{equation}
The sharp Hardy coefficient on $\Om^h$ is $C(N)$ for every $h>0$.
\end{lemma}
\begin{proof}
The definition gives
$\norm{\mathcal U_hu}_2^2=h^N\sum_n|u(hn)|^2=\norm{u}_{2,h}^2$,
and the inverse map is
$(\mathcal U_h^{-1}v)(hn)=h^{-N/2}v(n)$.
Substitution in \eqref{bg:mesh-laplacian} shows that the Laplacian acquires
the factor $h^{-2}$. The potential acquires the same factor because
$|hn|^{-2}=h^{-2}|n|^{-2}$. This proves
\eqref{bg:unitary-scaling}.

For completeness, the cancellation is also visible directly in the
variational quotient. If $v(n)=u(hn)$, define
\[
 \mathcal E_h[u]
 =h^{N-2}\sum_{n\in\Z^N}\sum_{j=1}^N
       |u(hn+he_j)-u(hn)|^2,
 \qquad
 M_h[u]=h^N\sum_{n\in\Om}\frac{|u(hn)|^2}{|hn|^2}.
\]
Then $\mathcal E_h[u]=h^{N-2}\E[v]$ and
$M_h[u]=h^{N-2}\M[v]$. The correspondence between finitely supported
$u$ and $v$ is bijective, so the infimum of
$\mathcal E_h[u]/M_h[u]$ equals $C(N)$.
\end{proof}

\subsection{Positivity and the bottom of the spectrum}

Let $\operatorname{Adj}_{\Om}$ denote the adjacency operator of the
punctured graph, so that $L_N=2NI-\operatorname{Adj}_{\Om}$.
The graph $\Om$ is connected. Indeed, a nearest-neighbor path in $\Z^N$
which meets the origin can be diverted around it through a second
coordinate direction. The following elementary bound will also be used
quantitatively. Fix $o=e_1$, and let $d_{\Om}$ be graph distance in $\Om$.
Then
\begin{equation}\label{bg:path-distance}
 d_{\Om}(n,o)\le |n|_1+3\qquad(n\in\Om).
\end{equation}
To verify this, start with a shortest path in the full lattice, whose
length is $|n-e_1|_1\le |n|_1+1$. If it contains a segment $a,0,b$, with
$a$ and $b$ distinct neighbors of the origin, replace this segment by
$a,a+b,b$ when $b\ne-a$. When $b=-a$, choose a coordinate unit vector
$e_j$ perpendicular to $a$ and use
$a,a+e_j,e_j,-a+e_j,-a$. The latter replacement adds two edges and every
new vertex is nonzero. A shortest path meets the origin at most once.

Connectedness ensures that the semigroup propagates positivity to every vertex.

\begin{lemma}\label{bg:positive-semigroup}
For every $\lambda\ge0$, $H_{\lambda,N}$ is bounded and selfadjoint on
$\ell^2(\Om)$, with
\[
 H_{\lambda,N}\ge-\lambda I,
 \qquad \norm{H_{\lambda,N}}\le4N+\lambda.
\]
The operators
$S_{\lambda,N}(s)=e^{-sH_{\lambda,N}}$ form a positivity improving
semigroup. For each initial value in $\ell^2(\Om)$, they give the unique
global continuously differentiable solution of the corresponding Cauchy
problem. The same assertions hold on every mesh $h$, with the bounds
scaled by $h^{-2}$.
\end{lemma}
\begin{proof}
The adjacency operator is symmetric and has norm at most $2N$. One way
to see the bound is to apply Cauchy--Schwarz at each vertex:
\[
 \sum_n\left|\sum_{m\sim n,\ m\ne0}v(m)\right|^2
 \le2N\sum_n\sum_{m\sim n,\ m\ne0}|v(m)|^2
 \le(2N)^2\norm v_2^2.
\]
The identity $\ip{v}{L_Nv}=\E[v]\ge0$, first obtained for finite
support, extends to $\ell^2$ by continuity. Also $0\le V\le I$ and
$\norm V=1$. These observations prove boundedness, selfadjointness and
the stated inequalities.

Set $B_\lambda=\operatorname{Adj}_{\Om}+\lambda V$. Its matrix entries are
nonnegative, and
\begin{equation}\label{bg:positive-power-series}
 S_{\lambda,N}(s)=e^{-2Ns}
       \sum_{k=0}^{\infty}\frac{s^k}{k!}B_\lambda^k
\end{equation}
converges in operator norm. Every term preserves nonnegativity. If
$n,m\in\Om$ and $d=d_{\Om}(n,m)$, a path of length $d$ gives
$(B_\lambda^d)(n,m)\ge1$. Thus the matrix entry of
$S_{\lambda,N}(s)$ from $m$ to $n$ is strictly positive for every $s>0$.
A nonzero nonnegative sequence has a strictly positive value at some
vertex, which proves positivity improvement.

The norm-convergent exponential of a bounded operator is differentiable
and solves the Cauchy problem. Uniqueness follows, for example, by
differentiating $e^{sH_{\lambda,N}}v(s)$ for any differentiable solution
$v$. Equation~\eqref{bg:unitary-scaling} transfers these assertions to
the mesh $h$.
\end{proof}

This also gives a direct proof of the positivity needed here. In the
general graph theory, the corresponding implication from connectedness
is \cite[Theorem~1.26]{KLW21}. Its Dirichlet-form hypotheses apply to
our attractive operator after the harmless shift
\[
 H_{\lambda,N}+\lambda I=L_N+\lambda(1-V),
\]
whose killing term is nonnegative.

The essential spectrum will ensure that a negative spectral bottom is an isolated eigenvalue.

\begin{lemma}\label{bg:essential-spectrum}
For every $\lambda\ge0$,
\begin{equation}\label{bg:essential-band}
 \sigma_{\mathrm{ess}}(H_{\lambda,N})=[0,4N].
\end{equation}
Consequently, $\inf\sigma(H_{\lambda,N})\le0$, and if this infimum is
negative, it is an isolated eigenvalue of finite multiplicity.
\end{lemma}
\begin{proof}
On $\ell^2(\Z^N)$, the Fourier transform identifies the full positive
Laplacian with multiplication by
\[
 \ell_N(\theta)=2N-2\sum_{j=1}^N\cos\theta_j,
 \qquad \theta\in[-\pi,\pi]^N.
\]
Its spectrum and essential spectrum are $[0,4N]$. Indeed, this is the
range of the continuous function $\ell_N$, and every neighborhood of
each value in this range has a preimage of positive measure.

Under
$\ell^2(\Z^N)=\ell^2(\Om)\oplus\operatorname{span}\{\delta_0\}$,
the operator $L_N\oplus0$ differs from the full Laplacian only through
the entries incident to the origin and its diagonal entry at the
origin. Their difference therefore has finite rank. The essential
spectrum is unchanged by this perturbation and by adjoining a
one-dimensional summand, so
$\sigma_{\mathrm{ess}}(L_N)=[0,4N]$.

Finally, $V$ is compact: its restriction to $\{|n|\le R\}$ has finite
rank, and the norm of the remaining diagonal operator is at most
$R^{-2}$. Invariance of the essential spectrum under compact
selfadjoint perturbations proves \eqref{bg:essential-band}. Any
spectral point below zero is therefore a discrete eigenvalue. Since
the spectrum is nonempty, compact and contains zero in its essential
part, the final assertions follow.
\end{proof}

The Hardy inequality now determines exactly when exponential growth begins.
The broader relation between lower spectral bounds and positive solutions
is developed in \cite[Theorem~4.2]{KLVW15}. Here the precise contractivity
threshold follows directly from the Hardy inequality and spectral calculus.

\begin{proposition}\label{bg:spectral-alternative}
The following statements hold.
\begin{enumerate}
\item If $0\le\lambda\le C(N)$, then
$\inf\sigma(H_{\lambda,N})=0$ and
\begin{equation}\label{bg:contractive-norm}
 \norm{S_{\lambda,N}(s)}_{2\to2}=1\qquad(s\ge0).
\end{equation}
\item If $\lambda>C(N)$, then
\begin{equation}\label{bg:mu-definition}
 \inf\sigma(H_{\lambda,N})=-\mu_N(\lambda)<0
\end{equation}
is a simple isolated eigenvalue with a normalized eigenfunction
$\phi_{\lambda,N}(n)>0$ at every vertex. Moreover,
\begin{equation}\label{bg:supercritical-norm}
 \norm{S_{\lambda,N}(s)}_{2\to2}
       =e^{\mu_N(\lambda)s}\qquad(s\ge0).
\end{equation}
\end{enumerate}
In particular, contractivity for any one positive time, or for all
positive times, is equivalent to $\lambda\le C(N)$.
\end{proposition}
\begin{proof}
The Hardy inequality and continuity of its two quadratic forms give
\[
 \ip{v}{H_{\lambda,N}v}
       =\E[v]-\lambda\M[v]\ge0
       \qquad(v\in\ell^2(\Om))
\]
whenever $\lambda\le C(N)$. Combining this with
Lemma~\ref{bg:essential-spectrum} gives the spectral bottom zero.

If $\lambda>C(N)$, the definition of the infimum supplies a nonzero
$v\in c_{00}(\Om)$ for which $\E[v]<\lambda\M[v]$. Its Rayleigh
quotient for $H_{\lambda,N}$ is negative. Therefore the spectral
bottom is negative and, by Lemma~\ref{bg:essential-spectrum}, is an
isolated eigenvalue of finite multiplicity.

We give the positivity and simplicity argument explicitly. Let $\psi$
be a normalized eigenfunction at this bottom, denoted by $E_0$.
Replacing $\psi$ by $|\psi|$ does not increase the energy, because
$\bigl||\psi(n)|-|\psi(m)|\bigr|\le|\psi(n)-\psi(m)|$ along every
edge, and it leaves the potential term and the norm unchanged.
Minimality of $E_0$ forces equality. Since $H_{\lambda,N}-E_0I$ is
nonnegative, equality of the quadratic form implies
$(H_{\lambda,N}-E_0I)|\psi|=0$. Thus there is a nonnegative
eigenfunction $\phi=|\psi|$. Its vertex equation is
\begin{equation}\label{bg:ground-state-equation}
 (2N+\mu_N(\lambda)-\lambda|n|^{-2})\phi(n)
       =\sum_{m\sim n,\ m\ne0}\phi(m).
\end{equation}
If $\phi$ vanished at one vertex, this equation would force it to
vanish at every neighbor and hence, by connectedness, everywhere.
Therefore $\phi>0$.

For any other eigenfunction at $E_0$, the preceding energy comparison
is again an equality. The difference of the two energies is the sum
over edges in $\Om$ of the nonnegative quantities
\[
 2\bigl(|\psi(n)||\psi(m)|
                 -\operatorname{Re}(\psi(n)\overline{\psi(m)})\bigr).
\]
The sum is finite or absolutely convergent, since $\psi\in\ell^2$ and
the vertex degrees are bounded. Every summand must vanish. The
absolute value of this eigenfunction is strictly positive by the
same argument as above, so its complex phase is constant along
edges, and consequently throughout $\Om$. Such an eigenfunction
cannot be orthogonal to $\phi$. An eigenspace of dimension greater
than one would contain a nonzero vector orthogonal to $\phi$, a
contradiction. This proves simplicity; compare
\cite[Lemma~5.2]{KLW21}.

For a bounded selfadjoint operator, spectral calculus gives
\[
 \norm{e^{-sH_{\lambda,N}}}
       =\sup_{\xi\in\sigma(H_{\lambda,N})}e^{-s\xi}
       =e^{-s\inf\sigma(H_{\lambda,N})}.
\]
This proves both norm identities and the final equivalence.
\end{proof}

Combining the preceding proposition with scaling yields the exact
mesh-dependent norm
\begin{equation}\label{bg:mesh-norm}
 \norm{e^{-tH_{\lambda,N}^h}}_{2,h\to2,h}
 =\begin{cases}
  1,&0\le\lambda\le C(N),\\[2pt]
  \exp\!\bigl(\mu_N(\lambda)t/h^2\bigr),&\lambda>C(N).
 \end{cases}
\end{equation}
Thus the stable assertion includes the critical value. It asserts
uniform contractivity, without requiring a classification of a
possible zero eigenfunction at that value.

At fixed mesh, the spectral theorem also describes the strong limit of the evolution.

\begin{remark}\label{bg:critical-limit}
For $0\le\lambda<C(N)$, $S_{\lambda,N}(s)v\to0$ in $\ell^2$ for every
$v\in\ell^2$ as $s\to\infty$. At $\lambda=C(N)$ the limit is the
orthogonal projection of $v$ onto $\ker H_{C(N),N}$. To justify both
claims, spectral calculus for a nonnegative selfadjoint operator gives
the strong limit $e^{-sH}\to\mathbf1_{\{0\}}(H)$ by dominated
convergence against each vector's spectral measure. In the strictly
subcritical case,
\[
 H_{\lambda,N}\ge\left(1-\frac{\lambda}{C(N)}\right)L_N.
\]
Hence a vector in its kernel has zero $L_N$-energy. Its zero extension
has equal values across every lattice edge and is zero at the origin,
so the vector is zero. This observation also explains why the norm
identity \eqref{bg:contractive-norm} does not preclude strong decay.
\end{remark}

\subsection{Diagonal growth and propagation along paths}

We next turn the norm growth into pointwise divergence for a fixed
macroscopic initial profile. The relevant strategy is the diagonal
growth and propagation argument in
\cite[Theorem~3.1]{GGKR24}. On the lattice, a single connecting path
provides the propagation estimate. This avoids both a parabolic
Harnack estimate and a bound on the spectral remainder.

Assume from now on that $\lambda>C(N)$, and abbreviate
$\mu=\mu_N(\lambda)$ and $\phi=\phi_{\lambda,N}$. Let
\[
 K_s(n,m)=\ip{\delta_n}{S_{\lambda,N}(s)\delta_m}
       \qquad(n,m\in\Om, s\ge0)
\]
be the semigroup kernel with respect to counting measure. It is
symmetric and nonnegative, and it is strictly positive for $s>0$.

Two kernel bounds will quantify growth at one vertex and propagation along a path.

\begin{lemma}\label{bg:kernel-bounds}
Fix $o=e_1$. For every $r\ge0$,
\begin{equation}\label{bg:diagonal-growth}
 K_r(o,o)\ge \phi(o)^2e^{\mu r}.
\end{equation}
For $n\in\Om$ and $\rho\ge d_{\Om}(n,o)$,
\begin{equation}\label{bg:path-kernel}
 K_\rho(n,o)\ge e^{-2N\rho}.
\end{equation}
Consequently, if
$\rho\ge\max\{d_{\Om}(n,o),d_{\Om}(m,o)\}$ and $s\ge2\rho$, then
\begin{equation}\label{bg:three-step-kernel}
 K_s(n,m)\ge
 \phi(o)^2\exp\!\bigl(\mu s-(4N+2\mu)\rho\bigr).
\end{equation}
\end{lemma}
\begin{proof}
Let $\nu_o$ be the positive spectral measure of $H_{\lambda,N}$
associated with $\delta_o$. The normalized eigenfunction $\phi$ gives
an atom of mass $\phi(o)^2$ at $-\mu$. All other contributions in
\[
 K_r(o,o)=\int_{\sigma(H_{\lambda,N})}e^{-r\xi}\,d\nu_o(\xi)
\]
are nonnegative. Keeping this atom proves
\eqref{bg:diagonal-growth}.

Put $d=d_{\Om}(n,o)$. The term of order $d$ in
\eqref{bg:positive-power-series} gives
\[
 K_\rho(n,o)\ge
 e^{-2N\rho}\frac{\rho^d}{d!}(B_\lambda^d)(n,o)
 \ge e^{-2N\rho}\frac{\rho^d}{d!}.
\]
For $d\ge1$ and $\rho\ge d$, one has
$\rho^d\ge d^d\ge d!$. For $d=0$, the term of order zero gives the
same assertion. This proves \eqref{bg:path-kernel}.

The semigroup identity and nonnegativity allow us to retain only the
intermediate vertex $o$ at both joining times:
\[
 K_s(n,m)
 \ge K_\rho(n,o)K_{s-2\rho}(o,o)K_\rho(o,m).
\]
These kernel compositions are legitimate because rows and columns of
a bounded operator belong to $\ell^2$; the relevant sums converge by
Cauchy--Schwarz. They can also be obtained first over finite sums and
then by monotone convergence. Applying the first two bounds yields
\eqref{bg:three-step-kernel}.
\end{proof}

\subsection{Fixed initial profiles and quantitative complete blow-up}

We combine these bounds at the time scale $t/h^2$ to obtain a rate uniform over bounded physical regions. This estimate adds the simultaneous dependence on $h$ to the fixed-vertex logarithmic asymptotics of \cite[Theorem~3.1]{KLVW15}.

\begin{proposition}\label{bg:pointwise-rate}
Let $\lambda>C(N)$ and let
$f\in C_c(\R^N\setminus\{0\})$ be nonnegative and not identically
zero. Give \eqref{bg:mesh-evolution} the sampled initial values
$u_{0,h}(hn)=f(hn)$. For every $R>0$ and $\tau>0$, there exist
$h_0>0$ and constants $c_R,c>0$ such that
\begin{equation}\label{bg:logarithmic-estimate}
 \mu_N(\lambda)t-c_Rh-ch^2
 \le h^2\log u_h(t,hn)
 \le \mu_N(\lambda)t+\frac N2h^2|\log h|+ch^2
\end{equation}
whenever $0<h\le h_0$, $t\ge\tau$ and $0<|hn|\le R$.
In particular,
\begin{equation}\label{bg:uniform-log-limit}
 \lim_{h\downarrow0}
 \sup_{\substack{t\ge\tau\\n\in\Om,\ |hn|\le R}}
 \left|h^2\log u_h(t,hn)-\mu_N(\lambda)t\right|=0,
\end{equation}
and
\begin{equation}\label{bg:complete-blowup}
 \lim_{h\downarrow0}
 \inf_{\substack{t\ge\tau\\n\in\Om,\ |hn|\le R}}
 u_h(t,hn)=+\infty.
\end{equation}
\end{proposition}
\begin{proof}
Set $f_h(n)=f(hn)$, viewed for the moment as a sequence on $\Om$
with counting measure. The two factors $h^{N/2}$ from the unitary
map and its inverse cancel. Consequently, scaling gives the exact
pointwise formula
\begin{equation}\label{bg:sampled-solution}
 u_h(t,hn)=\bigl(S_{\lambda,N}(t/h^2)f_h\bigr)(n)
       =\sum_{m\in\Om}K_{t/h^2}(n,m)f(hm).
\end{equation}
For each fixed $h$, the sequence $f_h$ has finite support.

Choose $y_0\ne0$ with $f(y_0)>0$. By continuity, there are $r_0>0$
and $a>0$ such that the closed ball $\overline{B(y_0,r_0)}$ avoids the
origin and $f\ge a$ on this ball. For sufficiently small $h$, choose
a nearest lattice point $hm_h$ to $y_0$. Then
$|hm_h-y_0|\le\sqrt N h/2<r_0$, so $m_h\in\Om$ and
\begin{equation}\label{bg:positive-sample}
 f(hm_h)\ge a,\qquad |hm_h|\le |y_0|+r_0.
\end{equation}

For $h\le1$, the path estimate \eqref{bg:path-distance} implies
\[
 \max\{d_{\Om}(n,o),d_{\Om}(m_h,o)\}\le\frac{D_R}{h},
 \qquad 0<|hn|\le R,
\]
where one may take
\begin{equation}\label{bg:distance-constant}
 D_R=\sqrt N\max\{R,|y_0|+r_0\}+3.
\end{equation}
Put $\rho_h=\lceil D_R/h\rceil$. Decrease $h_0$ so that the positive
sample above exists, $h_0\le\min\{1,R\}$, and
$h_0\le\tau/(2(D_R+1))$. Then, for every $t\ge\tau$ and
$0<h\le h_0$,
\[
 \frac{t}{h^2}\ge2\rho_h,
 \qquad \rho_h\le\frac{D_R}{h}+1.
\]
Retaining the single term $m=m_h$ in
\eqref{bg:sampled-solution} and applying
Lemma~\ref{bg:kernel-bounds} gives
\begin{equation}\label{bg:solution-lower-bound}
 u_h(t,hn)\ge
 a\phi(o)^2
 \exp\!\left(\frac{\mu t}{h^2}
                   -(4N+2\mu)\rho_h\right).
\end{equation}
This also shows directly that the logarithm in the proposition is
well defined. Taking logarithms and multiplying by $h^2$ yields
\[
 h^2\log u_h(t,hn)
 \ge\mu t-(4N+2\mu)D_Rh
       -\bigl(4N+2\mu+|\log(a\phi(o)^2)|\bigr)h^2.
\]

For the upper bound, choose $M>0$ such that
$\supp f\subset[-M,M]^N$. Counting the lattice points in this cube
gives, for $0<h\le1$,
\[
 \norm{f_h}_2^2
 \le\norm f_\infty^2(2M/h+1)^N
 \le\norm f_\infty^2(2M+1)^Nh^{-N}.
\]
Write $B_f=\norm f_\infty(2M+1)^{N/2}>0$. Since
$\norm{\delta_n}_2=1$, the norm formula
\eqref{bg:supercritical-norm} and \eqref{bg:sampled-solution} imply
\begin{equation}\label{bg:solution-upper-bound}
 u_h(t,hn)\le
 \norm{S_{\lambda,N}(t/h^2)f_h}_2
 \le B_fh^{-N/2}e^{\mu t/h^2}.
\end{equation}
Taking logarithms proves the upper estimate in
\eqref{bg:logarithmic-estimate}, after enlarging $c$ to cover
$|\log B_f|$ as well as the constant in the lower estimate.

All constants are independent of $n$, $t$ and $h$ in the stated
range. The errors in \eqref{bg:logarithmic-estimate} tend uniformly
to zero, proving \eqref{bg:uniform-log-limit}. Finally,
\eqref{bg:solution-lower-bound} is bounded below, uniformly for
$t\ge\tau$ and $|hn|\le R$, by
\[
 a\phi(o)^2\exp\!\left(
       \frac{\mu\tau}{h^2}
       -\frac{(4N+2\mu)D_R}{h}-(4N+2\mu)\right).
\]
Its logarithm tends to $+\infty$, since $\mu\tau>0$ and the positive
term of order $h^{-2}$ dominates the negative term of order $h^{-1}$.
This proves \eqref{bg:complete-blowup}.
\end{proof}

Equation~\eqref{bg:mesh-norm} and
Proposition~\ref{bg:pointwise-rate} prove
Theorem~\ref{thm:main-bg}. In particular, for any $x\ne0$ and any
sequence $x_h\in\Om^h$ with $x_h\to x$, the logarithmic limit holds
at $(t,x_h)$ for every fixed $t>0$. The estimates are stronger: they
hold uniformly on bounded sets of observation points and for all
times bounded away from zero.

The approximation parameter in the classical result of
Baras--Goldstein is a truncation of the singular potential
\cite[Proposition~2.1 and Theorem~2.2]{BG84}. Here it is the mesh
size. For each fixed mesh, Lemma~\ref{bg:positive-semigroup} gives
a global solution for every coupling. The conclusion
\eqref{bg:complete-blowup} describes its singular behavior under
mesh refinement. Throughout, the critical coefficient is the
mesh-independent number $C(N)$ from Lemma~\ref{bg:scaling}.
\section{Convergence to the continuum heat equation}\label{cont:section}

Below the mesh-independent contractivity threshold, the lattice evolution
approximates the continuum evolution uniformly over the entire time
half-line. We compare functions on the lattice with
functions on $\R^N$ by constant extension over cells. The limiting
operator is the Friedrichs realization of the continuum Hardy operator.
The proof identifies its resolvent directly, using the strict discrete
Hardy inequality to control first differences. This avoids any
assumption of regularity at the singularity.

\subsection{The comparison maps and the continuum realization}

For $n\in\Z^N$ and $h>0$, let
\[
 Q_h(n)=h\bigl(n+[-\tfrac12,\tfrac12)^N\bigr).
\]
These cells partition $\R^N$ up to their boundaries. For
$u\in\ell_h^2(h\Om)$, retain the convention $u(0)=0$ and define
\begin{equation}\label{cont:comparison}
 \begin{split}
 (I_hu)(x)&=u(hn),\qquad x\in Q_h(n),\quad n\in\Z^N,\\
 (P_hf)(hn)&=h^{-N}\int_{Q_h(n)}f(x)\dd x,
                  \qquad n\in\Om.
 \end{split}
\end{equation}
Thus $I_hu$ vanishes on the cell containing the origin, and no average
over that cell is retained by $P_h$. Direct integration gives
\begin{equation}\label{cont:comparison-identities}
 \norm{I_hu}_2=\norm{u}_{2,h},\qquad
 P_h=I_h^*,\qquad P_hI_h=I,\qquad \norm{P_h}\le1.
\end{equation}
Consequently $\Pi_h=I_hP_h$ is the orthogonal projection onto the range
of $I_h$. Moreover,
\begin{equation}\label{cont:projection-limit}
 \Pi_hf\longrightarrow f\quad\text{in }L^2(\R^N)
 \qquad(f\in L^2(\R^N)).
\end{equation}
Indeed, for a compactly supported continuous function the averages on
all nonzero cells converge to the function in $L^2$, by uniform
continuity, while the integral of its square over $Q_h(0)$ tends to
zero. Density and $\norm{\Pi_h}\le1$ prove the assertion for every
$L^2$ function.

Write
\[
 \mathcal D=C_c^\infty(\R^N\setminus\{0\}),\qquad
 A_N=\frac{(N-2)^2}{4}.
\]
For $0\le\lambda<A_N$, define the sesquilinear form
\begin{equation}\label{cont:continuum-form}
 \mathfrak q_\lambda(u,v)
 =\int_{\R^N}\nabla u\cdot\overline{\nabla v}\dd x
   -\lambda\int_{\R^N}\frac{u\overline v}{|x|^2}\dd x,
 \qquad \Dom(\mathfrak q_\lambda)=H^1(\R^N),
\end{equation}
and put $\mathfrak q_\lambda[u]=\mathfrak q_\lambda(u,u)$.
The continuum Hardy inequality gives
\begin{equation}\label{cont:continuum-coercivity}
 \left(1-\frac{\lambda}{A_N}\right)\norm{\nabla u}_2^2
 \le\mathfrak q_\lambda[u]\le\norm{\nabla u}_2^2.
\end{equation}
Hence the form is closed and nonnegative. Its associated selfadjoint
operator is denoted by $\mathsf H_{\lambda,N}^{F}$, and
\[
 \mathsf S_\lambda(t)=e^{-t\mathsf H_{\lambda,N}^{F}}.
\]
Equivalently, $\mathsf H_{\lambda,N}^{F}$ is the Friedrichs extension
of $-\Delta-\lambda|x|^{-2}$ initially defined on $\mathcal D$.

We record why the puncture does not impose an extra boundary condition
on this subcritical form domain. The space $\mathcal D$ is dense in
$H^1(\R^N)$ when $N\ge3$. To prove this, first approximate by
$g\in C_c^\infty(\R^N)$. Choose smooth $\eta_\varepsilon$ which is
zero on $B_\varepsilon$, equals one outside $B_{2\varepsilon}$, and
satisfies $|\nabla\eta_\varepsilon|\le c/\varepsilon$. For fixed $g$,
\[
 \norm{(1-\eta_\varepsilon)g}_2^2\le c_g\varepsilon^N,
 \qquad
 \norm{\nabla((1-\eta_\varepsilon)g)}_2^2
       \le c_g(\varepsilon^N+\varepsilon^{N-2}).
\]
Both quantities tend to zero. The equivalence of the form norm with the
$H^1$ norm in \eqref{cont:continuum-coercivity} proves that $\mathcal D$
is a form core. In particular, this assertion does not require
$\mathcal D$ to be an operator core.

The comparison with $A_N$ needed below follows directly from smooth
test functions, independently of any classification by dimension.

\begin{lemma}\label{cont:threshold-upper}
For every $N\ge3$, one has $C(N)\le A_N$.
\end{lemma}
\begin{proof}
Fix $0\ne\varphi\in\mathcal D$, and sample it on the lattice:
$\varphi_h(hn)=\varphi(hn)$. For sufficiently small $h$, this is a
nonzero finitely supported function which vanishes near the origin.
Riemann sums and convergence of first difference quotients give
\[
 \mathcal E_h[\varphi_h]\longrightarrow
       \int_{\R^N}|\nabla\varphi|^2\dd x,
 \qquad
 M_h[\varphi_h]\longrightarrow
       \int_{\R^N}|x|^{-2}|\varphi|^2\dd x.
\]
For example, the first limit follows by writing the energy as the
Riemann sum of
$\sum_j|h^{-1}(\varphi(x+he_j)-\varphi(x))|^2$; these functions
converge uniformly to $|\nabla\varphi|^2$ and have support in a fixed
compact set. By Lemma~\ref{bg:scaling},
\[
 C(N)\le
 \frac{\int|\nabla\varphi|^2}{\int|x|^{-2}|\varphi|^2}.
\]
The infimum of the quotient on the right over $\mathcal D\setminus\{0\}$
is $A_N$. For completeness, its lower bound follows by expanding
\[
 0\le\int\left|\nabla\varphi+
       \frac{N-2}{2}\frac{x}{|x|^2}\varphi\right|^2
 =\int|\nabla\varphi|^2-A_N\int|x|^{-2}|\varphi|^2.
\]
To see sharpness, fix $0\ne\eta\in C_c^\infty(\R)$ and set
\[
 \varphi_L(x)=|x|^{-(N-2)/2}\eta\left(\frac{\log|x|}{L}\right),
 \qquad L>0.
\]
Each $\varphi_L$ belongs to $\mathcal D$. Radial integration, followed
by $s=\log r$, shows that its quotient is
\[
 A_N+\frac1{L^2}\frac{\int_\R|\eta'(s)|^2\dd s}
                               {\int_\R|\eta(s)|^2\dd s}.
\]
Letting $L$ tend to infinity proves the claim.
\end{proof}

\subsection{The subcritical convergence theorem}

\begin{theorem}\label{cont:subcritical}
Let $N\ge3$ and $0\le\lambda<C(N)$. For every $f\in L^2(\R^N)$,
\begin{equation}\label{cont:all-time-convergence}
 \lim_{h\downarrow0}\ \sup_{t\ge0}
 \norm{I_hS_{\lambda,h}(t)P_hf-\mathsf S_\lambda(t)f}_2=0.
\end{equation}
More generally, the same assertion holds with $P_hf$ replaced by any
family $u_{0,h}\in\ell_h^2(h\Om)$ satisfying $I_hu_{0,h}\to f$ in
$L^2(\R^N)$. In particular, it holds for
$u_{0,h}(hn)=f(hn)$ when $f\in C_c(\R^N\setminus\{0\})$.
\end{theorem}

The continuum limit is therefore the solution
$u(t)=\mathsf S_\lambda(t)f$ of
\[
 \partial_tu=\Delta u+\lambda|x|^{-2}u,\qquad u(0)=f,
\]
understood through the Friedrichs realization. The convergence is
strong convergence for each fixed initial datum; it does not assert
convergence in operator norm. We give the proof in successive steps.

\begin{lemma}\label{cont:consistency}
For $\varphi\in\mathcal D$, let
$\varphi_h(hn)=\varphi(hn)$. For every fixed $\lambda\ge0$,
\begin{equation}\label{cont:core-consistency}
 \begin{split}
 I_h\varphi_h&\longrightarrow\varphi,\\
 I_hH_{\lambda,N}^{h}\varphi_h
    &\longrightarrow(-\Delta-\lambda|x|^{-2})\varphi
 \end{split}
 \qquad\text{in }L^2(\R^N).
\end{equation}
\end{lemma}
\begin{proof}
The first limit is the usual convergence of sampled continuous functions
under constant extension. The support of $\varphi$ is compact and has
positive distance from the origin. For sufficiently small $h$, both
$H_{\lambda,N}^{h}\varphi_h$ and the sampled continuum expression
are supported in one fixed compact set away from the origin. On that
set, Taylor's formula gives, uniformly at lattice points,
\[
 -\Delta_h\varphi_h(hn)=-\Delta\varphi(hn)+O(h^2).
\]
The potential term is exactly
$-\lambda|hn|^{-2}\varphi(hn)$. Near the origin all these functions
vanish, including the discrete expressions at vertices adjacent to the
origin. The uniform Taylor error has $\ell_h^2$ norm tending to zero,
and constant extension of the smooth continuum expression converges
in $L^2$. This proves the second limit.
\end{proof}

\begin{lemma}\label{cont:resolvent}
Under the assumptions of Theorem~\ref{cont:subcritical}, for every
$f\in L^2(\R^N)$,
\begin{equation}\label{cont:resolvent-convergence}
 I_h(H_{\lambda,N}^{h}+i)^{-1}P_hf
       \longrightarrow(\mathsf H_{\lambda,N}^{F}+i)^{-1}f
       \quad\text{in }L^2(\R^N).
\end{equation}
The same statement holds with $i$ replaced by $-i$.
\end{lemma}
\begin{proof}
Set
\[
 u_h=(H_{\lambda,N}^{h}+i)^{-1}P_hf,
 \qquad w_h=I_hu_h,\qquad
 \kappa=1-\frac{\lambda}{C(N)}>0.
\]
The mesh quadratic form is
\[
 \mathfrak q_{\lambda,h}[v]
 =\ip{H_{\lambda,N}^{h}v}{v}_{2,h}
 =\mathcal E_h[v]-\lambda M_h[v].
\]
The Hardy inequality, initially obtained on finite support, extends to
all of $\ell_h^2$ because the operators involved are bounded for each
fixed $h$. It yields
\begin{equation}\label{cont:mesh-coercivity}
 \mathfrak q_{\lambda,h}[v]\ge\kappa\mathcal E_h[v].
\end{equation}
Since $H_{\lambda,N}^{h}$ is nonnegative and selfadjoint, its resolvent
at $-i$ has norm at most one. Taking the real and imaginary parts of
the resolvent equation paired with $u_h$ gives
\begin{align}
 \norm{u_h}_{2,h}&\le\norm f_2,\label{cont:resolvent-bound}\\
 \mathfrak q_{\lambda,h}[u_h]
    &=\operatorname{Re}\ip{P_hf}{u_h}_{2,h}\le\norm f_2^2,
       \label{cont:energy-bound}\\
 \norm{u_h}_{2,h}^2
    &=\operatorname{Im}\ip{P_hf}{u_h}_{2,h}.
       \label{cont:imaginary-identity}
\end{align}
In particular, $\mathcal E_h[u_h]\le\kappa^{-1}\norm f_2^2$.

We next identify the weak derivatives of every weak limit of $w_h$.
For a function on $\R^N$, write
\[
 D_{h,j}w(x)=\frac{w(x+he_j)-w(x)}h,
 \qquad
 D_{-h,j}\psi(x)=\frac{\psi(x)-\psi(x-he_j)}h.
\]
The cell geometry and the extension by zero at the origin give the exact
identity
\begin{equation}\label{cont:difference-energy}
 \sum_{j=1}^N\norm{D_{h,j}w_h}_2^2=\mathcal E_h[u_h].
\end{equation}
Here the sum defining the energy includes edges incident to the
origin, precisely as in \eqref{eq:energy-mass}.

Consider any sequence of mesh sizes tending to zero. Boundedness in
$L^2(\R^N)^{N+1}$ permits extraction of a subsequence for which
$w_h\rightharpoonup u$ and $D_{h,j}w_h\rightharpoonup g_j$ in $L^2$.
For every $\psi\in C_c^\infty(\R^N)$, translation of the integral gives
\[
 \int D_{h,j}w_h\,\overline\psi
       =-\int w_h\,\overline{D_{-h,j}\psi}.
\]
Since $D_{-h,j}\psi\to\partial_j\psi$ strongly in $L^2$, passage to
the limit proves $g_j=\partial_j u$ in distributions. Thus
\begin{equation}\label{cont:limit-H1}
 u\in H^1(\R^N),\qquad
 \norm{\nabla u}_2^2\le\liminf_{h\to0}\mathcal E_h[u_h].
\end{equation}
This step uses first differences of the constant extensions; the
extensions themselves need not belong to $H^1$.

Fix $\varphi\in\mathcal D$ and use its samples $\varphi_h$ as tests
in the resolvent equation. Selfadjointness gives
\[
 \ip{u_h}{H_{\lambda,N}^{h}\varphi_h}_{2,h}
      +i\ip{u_h}{\varphi_h}_{2,h}
       =\ip{P_hf}{\varphi_h}_{2,h}.
\]
Lemma~\ref{cont:consistency} and weak convergence imply
\begin{equation}\label{cont:weak-resolvent}
 \mathfrak q_\lambda(u,\varphi)+i\ip u\varphi
       =\ip f\varphi\qquad(\varphi\in\mathcal D).
\end{equation}
Indeed, integration by parts is legitimate for these tests, and the
potential is bounded on their support. Since $u\in H^1$ and the Hardy
inequality bounds the potential pairing by
$A_N^{-1}\norm{\nabla u}_2\norm{\nabla\varphi}_2$, both sides of
\eqref{cont:weak-resolvent} are continuous in the $H^1$ norm of
$\varphi$. The form-core property extends this equation to all
$\varphi\in H^1(\R^N)$. By the definition of the operator associated
with the closed form,
\[
 u\in\Dom(\mathsf H_{\lambda,N}^{F}),\qquad
 (\mathsf H_{\lambda,N}^{F}+i)u=f.
\]
The solution is unique. Every weakly convergent subsequence therefore
has this same limit, proving weak convergence of the entire family.

Finally, \eqref{cont:imaginary-identity} and the adjoint relation
$P_h=I_h^*$ imply
\[
 \norm{w_h}_2^2
   =\operatorname{Im}\ip f{w_h}
   \longrightarrow\operatorname{Im}\ip fu
   =\norm u_2^2.
\]
Weak convergence together with convergence of norms proves strong
convergence. Replacing $i$ by $-i$ reverses the signs of both
imaginary-part identities and gives the other assertion.
\end{proof}

We include the passage from resolvents to heat evolution to make clear
that convergence is uniform also at $t=0$. This is an instance of the
general relation between resolvent and semigroup convergence, including
the framework of varying Hilbert spaces in \cite[Theorem~2.4]{KS03}.

\begin{lemma}\label{cont:resolvent-to-heat}
Let $H_h\ge0$ be selfadjoint operators on $\ell_h^2$, and let $H\ge0$
be selfadjoint on $L^2(\R^N)$. Suppose that
\[
 I_h(H_h\mathbin{\pm}i)^{-1}P_hf
        \longrightarrow(H\mathbin{\pm}i)^{-1}f
 \qquad(f\in L^2).
\]
Then, for every $T<\infty$ and $f\in L^2$,
\begin{equation}\label{cont:compact-time}
 \sup_{0\le t\le T}
 \norm{I_he^{-tH_h}P_hf-e^{-tH}f}_2\longrightarrow0.
\end{equation}
If additionally $\ker H=\{0\}$, the supremum in
\eqref{cont:compact-time} can be taken over all $t\ge0$.
\end{lemma}
\begin{proof}
Identify the lattice space with the closed range of $I_h$, and extend
$H_h$ by zero on its orthogonal complement. The resulting selfadjoint
operator on $L^2(\R^N)$ is denoted by $\widetilde H_h$. In the present
application it is simply the bounded operator $I_hH_hP_h$. Its
resolvents and heat semigroup satisfy
\begin{align}
 (\widetilde H_h+i)^{-1}
    &=I_h(H_h+i)^{-1}P_h-i(I-\Pi_h),\label{cont:lifted-resolvent}\\
 e^{-t\widetilde H_h}
    &=I_he^{-tH_h}P_h+(I-\Pi_h).\label{cont:lifted-semigroup}
\end{align}
There is an analogous resolvent identity with $-i$. By
\eqref{cont:projection-limit}, the lifted resolvents converge strongly
to those of $H$.

Here is a direct functional-calculus argument for the required heat
convergence. On $[0,\infty)$ the algebra generated by
$(s+i)^{-1}$ and $(s-i)^{-1}$, without a constant term, is uniformly
dense in $C_0([0,\infty))$ by the Stone--Weierstrass theorem. Products
of the corresponding resolvents converge strongly, because all
factors are uniformly bounded. Uniform approximation and the spectral
theorem therefore give
\[
 F(\widetilde H_h)f\longrightarrow F(H)f
 \qquad(F\in C_0([0,\infty)),\ f\in L^2).
\]
Taking $F(s)=e^{-ts}$ proves convergence at each fixed $t>0$; at $t=0$
the lifted semigroups are the identity.

To obtain uniformity on bounded time intervals, first use data
$g=(H+i)^{-1}k$ and let
$g_h=(\widetilde H_h+i)^{-1}k$. Then $g_h\to g$, and
\[
 \norm{\widetilde H_hg_h}_2
       =\norm{k-ig_h}_2\le2\norm k_2.
\]
For $s,t\ge0$, integration of the semigroup derivative gives
\[
 \norm{e^{-t\widetilde H_h}g_h-e^{-s\widetilde H_h}g_h}_2
       \le2\norm k_2|t-s|.
\]
The analogous estimate holds for $e^{-tH}g$. Convergence at a finite
set of times and these bounds give uniform convergence on $[0,T]$
for the data $g_h,g$. The range of $(H+i)^{-1}$ is dense in $L^2$.
Approximation by such $g$, convergence $g_h\to g$, and contractivity
then prove uniform convergence for any fixed datum $f$. Formula
\eqref{cont:lifted-semigroup} and \eqref{cont:projection-limit} prove
\eqref{cont:compact-time}.

If $\ker H=\{0\}$, the spectral theorem gives
$e^{-tH}f\to0$ as $t\to\infty$. For $t\ge T$, contractivity yields
\[
 \norm{I_he^{-tH_h}P_hf-e^{-tH}f}_2
 \le\norm{e^{-TH_h}P_hf}_{2,h}+\norm{e^{-TH}f}_2.
\]
Convergence at time $T$ consequently implies
\[
 \limsup_{h\downarrow0}\ \sup_{t\ge T}
 \norm{I_he^{-tH_h}P_hf-e^{-tH}f}_2
       \le2\norm{e^{-TH}f}_2.
\]
Let $T$ tend to infinity, and combine this estimate with convergence on
$[0,T]$.
\end{proof}

\begin{proof}[Proof of Theorem~\ref{cont:subcritical}]
Lemma~\ref{cont:threshold-upper} gives $\lambda<A_N$. The continuum
operator is therefore the operator associated with
\eqref{cont:continuum-form}. Its kernel is zero: if
$\mathsf H_{\lambda,N}^{F}u=0$, then
$\mathfrak q_\lambda[u]=0$, and
\eqref{cont:continuum-coercivity} gives $\nabla u=0$. The only constant
function in $L^2(\R^N)$ is zero. Lemmas~\ref{cont:resolvent} and
\ref{cont:resolvent-to-heat} prove \eqref{cont:all-time-convergence}.

For a more general family of initial data, the isometry of $I_h$ gives
\[
 \norm{u_{0,h}-P_hf}_{2,h}
 \le\norm{I_hu_{0,h}-f}_2+\norm{f-\Pi_hf}_2\longrightarrow0.
\]
The lattice semigroups are contractions, so this difference remains
bounded by the same quantity for every $t\ge0$. This proves the
second assertion. If $f\in C_c(\R^N\setminus\{0\})$, its constant
extensions from sampled values converge to $f$ in $L^2$, by uniform
continuity on a fixed compact set, proving the final assertion.
\end{proof}

\subsection{An endpoint consequence of essential selfadjointness}

At $\lambda=C(N)$ the coefficient $\kappa$ in
\eqref{cont:mesh-coercivity} vanishes. Nonnegativity of the mesh
operators alone does not supply the $H^1$ bound used in
Lemma~\ref{cont:resolvent}. Nevertheless, when the continuum minimal
operator is essentially selfadjoint, consistency on its operator core
identifies the limit without that estimate.

\begin{proposition}\label{cont:endpoint-esa}
Suppose $N\ge3$ and
\begin{equation}\label{cont:endpoint-condition}
 C(N)\le A_N-1.
\end{equation}
Then all conclusions of Theorem~\ref{cont:subcritical} hold also for
$\lambda=C(N)$. In particular, the endpoint is included for every
$N\ge10$.
\end{proposition}
\begin{proof}
The classical essential-selfadjointness criterion for the inverse-square
operator states that $-\Delta-\lambda|x|^{-2}$ on $\mathcal D$ is
essentially selfadjoint when $\lambda\le A_N-1$; see
\cite[Theorem~8.1]{FMT07} and the original results cited there. One
way to see the threshold is to separate spherical harmonics. For a
spherical harmonic of degree $\ell$, the unitary radial transformation
gives the one-dimensional expression
\[
 -\frac{d^2}{dr^2}+\frac{\nu_\ell^2-1/4}{r^2},
 \qquad
 \nu_\ell^2=A_N+\ell(\ell+N-2)-\lambda.
\]
When $\lambda\le A_N-1$, every $\nu_\ell\ge1$. The zero-energy
solutions at the origin have powers $r^{1/2+\nu_\ell}$ and
$r^{1/2-\nu_\ell}$, and the latter is not square integrable there.
The origin is therefore limit point in every channel, as is infinity.
The resulting unique selfadjoint realization agrees with the Friedrichs
realization. In particular, $\mathcal D$ is an operator core for
$\mathsf H_{\lambda,N}^{F}$.

We give the consistency argument explicitly. Put $\lambda=C(N)$ and
$H=\mathsf H_{\lambda,N}^{F}$. Since $\mathcal D$ is an operator
core, the set
\[
 \{(H+i)\varphi:\varphi\in\mathcal D\}
\]
is dense in $L^2$: approximate any vector in $\Dom(H)$ in the graph
norm and use the surjectivity of $H+i$. For
$f=(H+i)\varphi$ in this dense set, let $\varphi_h$ be the samples
of $\varphi$. The mesh operators are nonnegative also at the endpoint,
so their resolvents at $-i$ have norm at most one. Hence
\begin{align*}
 &\norm{I_h(H_{\lambda,N}^{h}+i)^{-1}P_hf-\varphi}_2\\
 &\quad\le
 \norm{P_hf-(H_{\lambda,N}^{h}+i)\varphi_h}_{2,h}
       +\norm{I_h\varphi_h-\varphi}_2\longrightarrow0.
\end{align*}
For the first term, apply the isometry $I_h$, use $\Pi_hf\to f$,
and then Lemma~\ref{cont:consistency}. Uniform boundedness of the
resolvents extends convergence to every $f\in L^2$. The same argument
with $-i$ proves the other resolvent convergence. Lemma
\ref{cont:resolvent-to-heat} now applies.

Condition~\eqref{cont:endpoint-condition} gives $\lambda<A_N$, so
\eqref{cont:continuum-coercivity} again shows that $\ker H=\{0\}$.
We obtain uniform convergence for all $t\ge0$, and the proof for general
and sampled data is unchanged. Finally,
Proposition~\ref{prop:coercivity} gives $C(N)\le2N$, while
\[
 A_N-1-2N=\frac{N(N-12)}4\ge0\qquad(N\ge12).
\]
Thus every dimension $N\ge12$ satisfies
\eqref{cont:endpoint-condition}. Lemma~\ref{cont:endpoint-dimensions}
below proves the same condition for $N=10,11$.
\end{proof}

Two finite tests make the endpoint criterion explicit in the remaining
small dimensions. The dimension-eleven test is a specialization of the
two-orbit variational method in \cite{LizamaLowDimensions}; the following
integer amplitudes provide direct rational bounds with no spectral
computation.

\begin{lemma}
\label{cont:endpoint-dimensions}
One has
\[
 C(10)\le\frac{576104000}{38613219}<15=A_{10}-1,
 \qquad C(11)\le\frac{772}{41}<\frac{77}{4}=A_{11}-1.
\]
Consequently $C(N)\le A_N-1$ for every $N\ge10$.
\end{lemma}

\begin{proof}
An orbit consists of all signed permutations of its indicated nonzero
coordinates. In dimension ten, write $P_k=(1^k)$ and $Q_j=(2,1^j)$,
with cardinalities
$p_k=2^k\binom{10}{k}$ and $q_j=2^{j+1}10\binom9j$.
Define an invariant test $u$ by the following integer amplitudes, and
set it equal to zero on every other orbit, including the origin:
\begin{center}
\small
\begin{tabular}{rrrrrrr}
\toprule
$k$ & $p_k$ & $u(P_k)=a_k$ &\quad& $j$ & $q_j$ & $u(Q_j)=b_j$\\
\midrule
1 & 20 & 500 && 0 & 20 & 51\\
2 & 180 & 139 && 1 & 360 & 18\\
3 & 960 & 44 && 2 & 2880 & 7\\
4 & 3360 & 16 && 3 & 13440 & 3\\
5 & 8064 & 6 && 4 & 40320 & 1\\
6 & 13440 & 2 && & &\\
\bottomrule
\end{tabular}
\end{center}
The only edges with two nonzero endpoint values join $P_k$ to $P_{k+1}$,
$Q_j$ to $Q_{j+1}$, or $P_k$ to $Q_{k-1}$. Their respective counts are
$2(10-k)p_k$, $2(9-j)q_j$, and $kp_k$: these follow by inserting a
signed unit coordinate, or increasing one unit coordinate to two.
Thus, with each unoriented edge counted once,
\begin{align*}
 B:=\|u\|_2^2
 &=\sum_{k=1}^6p_ka_k^2+\sum_{j=0}^4q_jb_j^2=12011624,\\
 M_{10}[u]
 &=\sum_{k=1}^6\frac{p_ka_k^2}{k}
   +\sum_{j=0}^4\frac{q_jb_j^2}{j+4}
   =\frac{38613219}{5},\\
 F:=\sum_{\{x,y\}:x\sim y}u(x)u(y)
 &=\sum_{k=1}^5 2(10-k)p_ka_ka_{k+1}
   +\sum_{j=0}^3 2(9-j)q_jb_jb_{j+1}
   +\sum_{k=1}^5 kp_ka_kb_{k-1}\\
 &=56934240+2386800+3184800=62505840.
\end{align*}
Since the full lattice has degree twenty,
$\mathcal E_{10}[u]=20B-2F=115220800$; the diagonal term includes
edges leaving the support. Hence
$\mathcal E_{10}[u]-15M_{10}[u]=-618857<0$ and the first bound follows.

In dimension eleven, take $u=6$ on $(1)$ and $u=1$ on $(1,1)$,
and zero elsewhere. The orbit sizes are $22$ and $220$, and there are
$22\cdot20=440$ connecting edges. Therefore
\[
 \mathcal E_{11}[u]=22(22\cdot36+220)-2\cdot440\cdot6=16984,
 \qquad M_{11}[u]=22\cdot36+220/2=902,
\]
which proves the second bound. Finally, a test supported at one vertex
gives $C(N)\le2N$, while
$A_N-1-2N=N(N-12)/4\ge0$ for $N\ge12$.
\end{proof}
\subsection{An endpoint argument when the lattice and continuum thresholds agree}
\label{subsec:endpoint-notes}

The loss of uniform coercivity at the sharp lattice coupling does not by
itself prevent a continuum limit.  A useful substitute is uniqueness of the
nonnegative selfadjoint extension of the limiting minimal operator.  We
first give the precise functional-analytic argument; it also explains why
the same proof cannot be applied indiscriminately to every lattice-critical
coupling.

The following dual form is a standard construction of the least nonnegative extension. We include its proof to make the endpoint argument self-contained.

\begin{lemma}
\label{lem:endpoint-dual}
Let $S$ be a densely defined nonnegative symmetric operator on a Hilbert
space $\mathcal H$.  Put $q_S[v]=\langle Sv,v\rangle$ for $v\in \Dom(S)$ and
define, with the value $+\infty$ permitted,
\begin{equation}\label{eq:endpoint-dual-form}
 Q_S[u]=\sup_{v\in \Dom(S)}
 \bigl\{2\operatorname{Re}\langle u,Sv\rangle-q_S[v]\bigr\}.
\end{equation}
Then $Q_S$ is a densely defined closed nonnegative quadratic form whose
associated operator is a nonnegative selfadjoint extension of $S$.
Moreover, for every nonnegative selfadjoint extension $B$ of $S$,
\[
 \Dom(B^{1/2})\subset \Dom(Q_S),\qquad
 Q_S[u]\le\|B^{1/2}u\|^2\quad(u\in \Dom(B^{1/2})).
\]
In particular, if $S$ has only one nonnegative selfadjoint extension,
then $Q_S$ is the Friedrichs form of $S$.
\end{lemma}

\begin{proof}
The Cauchy--Schwarz inequality for the nonnegative form $q_S$ shows that
$q_S[v]=0$ implies $Sv=0$: indeed,
$|\langle Sv,w\rangle|^2\le q_S[v]q_S[w]=0$ for every $w\in \Dom(S)$,
and $\Dom(S)$ is dense.  Thus the quotient of $\Dom(S)$ by $\ker S$, equipped
with the inner product induced by $q_S$, is a pre-Hilbert space.  Denote
its completion by $\mathcal E$ and the equivalence class of $v$ by $[v]$.

On the dense subspace $\{[v]:v\in \Dom(S)\}$ of $\mathcal E$ define
\[
 J[v]=Sv\in\mathcal H.
\]
The preceding observation makes $J$ well-defined.  For $u\in \Dom(S)$,
symmetry gives
\[
 \langle u,J[v]\rangle_{\mathcal H}
   =\langle [u],[v]\rangle_{\mathcal E}.
\]
Consequently $u\in \Dom(J^*)$ and $J^*u=[u]$.  In particular, $J^*$ is
densely defined and $J$ is closable.  Let $\overline J$ be its closure.
The form
\[
 Q[u]=\|J^*u\|_{\mathcal E}^2,
 \qquad \Dom(Q)=\Dom(J^*),
\]
is densely defined, closed and nonnegative, because $J^*$ is a closed
operator.  Its associated selfadjoint operator is
$\overline J J^*$.  If $u\in \Dom(S)$, then $J^*u=[u]\in \Dom(J)$, so
$\overline J J^*u=Su$.  Hence this operator extends $S$.

For $u\in\mathcal H$, boundedness of
$[v]\mapsto\langle u,Sv\rangle$ in the $\mathcal E$-norm is exactly
the condition $u\in \Dom(J^*)$.  Optimizing over scalar multiples of $v$
and then using density in $\mathcal E$ gives
\[
 \sup_{v\in \Dom(S)}
 \{2\operatorname{Re}\langle u,Sv\rangle-\|[v]\|_{\mathcal E}^2\}
 =\|J^*u\|_{\mathcal E}^2,
\]
with both sides interpreted as $+\infty$ outside $\Dom(J^*)$.  This proves
\eqref{eq:endpoint-dual-form}.

Finally, if $B$ is a nonnegative selfadjoint extension and
$u\in \Dom(B^{1/2})$, then, for $v\in \Dom(S)\subset \Dom(B)$,
\[
 2\operatorname{Re}\langle u,Sv\rangle-q_S[v]
 =2\operatorname{Re}\langle B^{1/2}u,B^{1/2}v\rangle
       -\|B^{1/2}v\|^2
 \le\|B^{1/2}u\|^2.
\]
Taking the supremum proves the asserted comparison.  If the extension is
unique, $\overline J J^*$ must be the Friedrichs extension, proving the
last assertion.  This argument also identifies $Q_S$ as the
Krein--von Neumann form, but that terminology is not needed below.
\end{proof}

Write $A_N=(N-2)^2/4$ and let
\[
 S_{A_N}=-\Delta-A_N|x|^{-2},
 \qquad \Dom(S_{A_N})=C_c^\infty(\mathbb R^N\setminus\{0\}).
\]
Its Friedrichs extension is denoted by $\mathsf H_{A_N,N}^{F}$ and its
closed form by $\mathfrak q_{A_N}^{F}$.  This is the closure of the
nonnegative critical Hardy form, extended by $+\infty$ outside its domain. On the completed domain this closed form is defined by approximation, since the two integrals in its initial expression need not be finite separately.

We next apply the known classification of nonnegative Bessel realizations to the continuum critical operator.

\begin{lemma}
\label{lem:endpoint-critical-unique}
For every $N\ge3$, the operator $S_{A_N}$ has exactly one nonnegative
selfadjoint extension.  It follows that
\begin{equation}\label{eq:endpoint-critical-dual}
 \mathfrak q_{A_N}^{F}[u]
 =\sup_{\psi\in C_c^\infty(\mathbb R^N\setminus\{0\})}
 \left\{2\operatorname{Re}\langle u,S_{A_N}\psi\rangle
       -\langle S_{A_N}\psi,\psi\rangle\right\}.
\end{equation}
\end{lemma}

\begin{proof}
Expand in spherical harmonics of degree $\ell$, whose spherical
Laplacian eigenvalue is $\ell(\ell+N-2)$, and conjugate the radial part
by multiplication by $r^{(N-1)/2}$. On each radial channel the minimal
operator $S_\ell$ acts by
\begin{equation}\label{eq:endpoint-bessel-channel}
 -\frac{d^2}{dr^2}+\frac{\nu_\ell^2-1/4}{r^2},
 \qquad \nu_\ell^2=\ell(\ell+N-2),
 \qquad \Dom(S_\ell)=C_c^\infty(0,\infty).
\end{equation}
The cancellation leading to this formula is
$(N-1)(N-3)/4-A_N=-1/4$.
For $\ell\ge1$, one has $\nu_\ell^2\ge N-1\ge2$, so these channels
are essentially selfadjoint.  For $\ell=0$, one has $\nu_0=0$.
Bruneau, Derezi\'nski and Georgescu prove that the Friedrichs and
Krein extensions of this last minimal Bessel operator coincide,
and that every other selfadjoint extension has a negative eigenvalue
\cite[Proposition~4.21; see also Theorem~4.22]{BDG11}.

For completeness, this channel reduction applies to the entire minimal
operator, rather than only to radial functions.  Finite spherical-harmonic
sums with smooth radial coefficients supported in a compact annulus form
an operator core for its closure: angular spectral truncation converges
together with the operator on each such annulus, where all coefficients
are smooth and bounded.  The sum of the channels with $\ell\ge1$ is
selfadjoint after closure.  Its deficiency spaces are zero, so the
deficiency spaces of the full minimal operator are those of the single
radial channel.  A selfadjoint extension therefore consists of an
extension in that radial channel and the fixed selfadjoint operator on
its orthogonal complement.  Nonnegativity selects the unique radial
extension just described.  This proves uniqueness for $S_{A_N}$.
Equation~\eqref{eq:endpoint-critical-dual} now follows from
Lemma~\ref{lem:endpoint-dual}.
\end{proof}

We retain the comparison maps $I_h,P_h$ from \eqref{cont:comparison} and
write $\mathfrak q_{A_N,h}[u]=\ip{H_{A_N,N}^{h}u}{u}_{2,h}$.

These facts identify the continuum limit also at a common lattice and continuum threshold.
For clarity, Mosco convergence under $I_h$ means the following two conditions:
weak convergence $I_hu_h\rightharpoonup u$ implies the lower bound on form values,
and every $u$ admits a strongly convergent family attaining the upper bound.
The proof below verifies both conditions directly.

\begin{theorem}
\label{thm:endpoint-critical-mosco}
Assume $C(N)=A_N$.  Under the embeddings $I_h$, the forms
$\mathfrak q_{A_N,h}$ converge in the Mosco sense to
$\mathfrak q_{A_N}^{F}$.  Consequently, for every $a>0$ and
$f\in L^2(\mathbb R^N)$,
\begin{equation}\label{eq:endpoint-resolvent-limit}
 I_h(H_{A_N,N}^{h}+a)^{-1}P_hf
 \longrightarrow(\mathsf H_{A_N,N}^{F}+a)^{-1}f
 \quad\text{in }L^2(\mathbb R^N).
\end{equation}
The heat semigroups converge uniformly for all $t\ge0$:
\begin{equation}\label{eq:endpoint-semigroup-limit}
 \sup_{t\ge0}\norm{I_he^{-tH_{A_N,N}^{h}}P_hf
             -e^{-t\mathsf H_{A_N,N}^{F}}f}_2\longrightarrow0.
\end{equation}
In particular, the theorem applies in every dimension $3\le N\le8$
by \eqref{eq:known-transition}. The conclusion also holds for every family
of initial data whose constant extensions converge strongly to $f$.
\end{theorem}

\begin{proof}
The equality of thresholds implies $\mathfrak q_{A_N,h}\ge0$ for every $h>0$.
For $\psi\in C_c^\infty(\mathbb R^N\setminus\{0\})$, let
$\psi_h(hn)=\psi(hn)$.  Its support is finite and stays away from the
removed cell for sufficiently small $h$.  Taylor's formula for central
differences and ordinary Riemann sums give
\begin{equation}\label{eq:endpoint-consistency}
 \begin{aligned}
 I_h\psi_h&\to\psi,\qquad
 I_hH_{A_N,N}^{h}\psi_h\to S_{A_N}\psi
       &&\text{strongly in }L^2,\\
 \mathfrak q_{A_N,h}[\psi_h]&\to\langle S_{A_N}\psi,\psi\rangle.
 \end{aligned}
\end{equation}
Here all functions involved have supports in a fixed compact annulus,
enlarged by at most one lattice step.  On that annulus the potential is
smooth and bounded.  The consistency error in the lattice differential
expression is uniformly $O(h^2)$, while the piecewise constant
interpolation error tends to zero uniformly.  These observations justify
each asserted strong convergence without an estimate near the singularity.

We first prove the lower-limit condition.  Suppose
$I_hu_h\rightharpoonup u$ weakly in $L^2$.  Nonnegativity applied to
$u_h-\psi_h$ gives
\[
 \mathfrak q_{A_N,h}[u_h]
 \ge 2\operatorname{Re}\langle I_hu_h,I_hH_{A_N,N}^{h}\psi_h\rangle
       -\mathfrak q_{A_N,h}[\psi_h].
\]
Using \eqref{eq:endpoint-consistency} and weak convergence, then taking
the supremum over $\psi$, yields
\[
 \liminf_{h\downarrow0}\mathfrak q_{A_N,h}[u_h]
 \ge \mathfrak q_{A_N}^{F}[u]
\]
by \eqref{eq:endpoint-critical-dual}.  This remains valid if the right
side is infinite.

For the upper-limit condition, take $u\in \Dom(\mathfrak q_{A_N}^{F})$.
By the definition of the Friedrichs closure there are
$\psi_j\in C_c^\infty(\mathbb R^N\setminus\{0\})$ converging to $u$
in the norm
$(\|\cdot\|_2^2+\mathfrak q_{A_N}^{F}[\cdot])^{1/2}$.
For each fixed $j$, \eqref{eq:endpoint-consistency} applies to its
lattice samples.  Choose a decreasing sequence $h_j\downarrow0$ such
that, for $0<h<h_j$, both the $L^2$ sampling error and the absolute
error in the form value of $\psi_j$ are less than $1/j$.
For $h_{j+1}\le h<h_j$ use the samples of $\psi_j$.
The resulting $u_h$ satisfy
$I_hu_h\to u$ strongly and
$\mathfrak q_{A_N,h}[u_h]\to \mathfrak q_{A_N}^{F}[u]$.
For $u$ outside the limiting form domain, the upper-limit requirement is automatic: take $u_h=P_hu$ and use \eqref{cont:projection-limit}.
This proves Mosco convergence.

The resolvent conclusion can also be read directly from these two
conditions.  The vector $(H_{A_N,N}^{h}+a)^{-1}P_hf$ uniquely minimizes
\[
 \mathfrak q_{A_N,h}[v]+a\norm v_{2,h}^2
       -2\operatorname{Re}\langle f,I_hv\rangle.
\]
The minimizers are bounded by $a^{-1}\|f\|_2$.  The lower-limit
condition and recovery sequences show that every weak limit minimizes
the corresponding continuum functional; that minimizer is unique.
The identities
\[
 \mathfrak q_{A_N,h}[u_h]+a\norm{u_h}_{2,h}^2
   =\operatorname{Re}\langle f,I_hu_h\rangle
\]
and the lower-limit condition give, for the identified weak limit $u$,
\[
 a\limsup_{h\downarrow0}\norm{u_h}_{2,h}^2
 \le \operatorname{Re}\langle f,u\rangle-\mathfrak q_{A_N}^{F}[u]
 =a\norm u_2^2.
\]
Weak lower semicontinuity gives the opposite inequality for the limit
inferior. The norms therefore converge, and weak convergence becomes
strong convergence, proving \eqref{eq:endpoint-resolvent-limit}.

To pass from positive resolvents to the heat evolution, observe that each
function $x\mapsto(x\pm i)^{-1}$ on $[0,\infty)$ can be approximated
uniformly by polynomials in $(x+a)^{-1}$ with zero constant term. Indeed,
after the substitution $y=(x+a)^{-1}$ it extends continuously to
$[0,a^{-1}]$, with value zero at $y=0$. The Weierstrass approximation
theorem, $P_hI_h=I$, and the uniform resolvent bounds therefore give
convergence of the two complex resolvents. Lemma~\ref{cont:resolvent-to-heat}
applies once the continuum kernel is shown to be zero.

Here this last fact follows from the critical ground-state identity.
Put $p=(N-2)/2$. For $\psi\in\mathcal D$, integration by parts gives
\[
 \mathfrak q_{A_N}^{F}[\psi]
   =\int_{\R^N}|x|^{-2p}\left|\nabla\bigl(|x|^p\psi\bigr)\right|^2\dd x.
\]
If $\mathsf H_{A_N,N}^{F}u=0$, choose core functions $\psi_j$ converging
to $u$ in the critical form norm. Their transformed gradients tend to
zero in the displayed weighted space. On every compact annulus,
$|x|^p\psi_j\to|x|^pu$ in $L^2$, so distributional differentiation
gives $\nabla(|x|^pu)=0$. The punctured space is connected, and hence
$u=c|x|^{-p}$ there. This function is not in $L^2(\R^N)$ unless $c=0$.
Thus the kernel is zero. Lemma~\ref{cont:resolvent-to-heat} proves
\eqref{eq:endpoint-semigroup-limit}, and contractivity treats strongly
compatible initial data exactly as in Theorem~\ref{cont:subcritical}.
\end{proof}

\subsection{Failure of convergence between the two thresholds}

The known transition \eqref{eq:known-transition} has a macroscopic
consequence: a nonempty interval of couplings admits a contractive
continuum evolution while the sampled lattice solutions diverge on every
nonempty bounded open observation set. The following statement combines the
macroscopic blow-up estimate proved here with the continuum Friedrichs
construction and gives a dynamical meaning to the gap between the two
Hardy constants.

\begin{corollary}\label{cont:mismatch}
Let $N\ge9$ and $C(N)<\lambda\le A_N$. For every nonzero nonnegative
$f\in C_c(\R^N\setminus\{0\})$, the continuum solution
$u(t)=e^{-t\mathsf H_{\lambda,N}^{F}}f$ satisfies
$\norm{u(t)}_2\le\norm f_2$. In contrast, the lattice solutions with
$u_h(0,hn)=f(hn)$ satisfy, for every nonempty bounded open set
$U\subset\R^N$ and every $\tau>0$,
\begin{equation}\label{cont:local-divergence}
 \inf_{t\ge\tau}\norm{I_hu_h(t)}_{L^2(U)}\longrightarrow+\infty
 \qquad(h\downarrow0).
\end{equation}
In particular, these solutions cannot converge locally in $L^2$ to the
continuum Friedrichs evolution at any positive time.
\end{corollary}
\begin{proof}
For $\lambda<A_N$ the continuum form is nonnegative by
\eqref{cont:continuum-coercivity}; for $\lambda=A_N$ use its Friedrichs
closure defined above. Both generate contraction semigroups.
Choose $R$ so that $U\subset B_{R-1}$. For small $h$, every cell meeting
$U$ has center of norm at most $R$. The quantitative lower bound in
Proposition~\ref{bg:pointwise-rate} gives constants $a_R,b$ independent
of $h,t$ such that, on the union of its nonzero cells,
\[
 I_hu_h(t,x)\ge
 \exp\bigl(\mu_N(\lambda)\tau/h^2-a_R/h-b\bigr)
 \qquad(t\ge\tau).
\]
The removed cell has measure $h^N$, so integration over $U$ gives
\[
 \norm{I_hu_h(t)}_{L^2(U)}^2
 \ge (|U|-h^N)
       \exp\bigl(2\mu_N(\lambda)\tau/h^2-2a_R/h-2b\bigr)
\]
for all sufficiently small $h$. Since $|U|>0$ and
$\mu_N(\lambda)>0$, this proves \eqref{cont:local-divergence}.
\end{proof}

At the lattice-critical coupling the preceding arguments establish
convergence whenever $C(N)=A_N$ or $C(N)\le A_N-1$. They do not identify
the limit in the intermediate case
\begin{equation}\label{cont:unresolved-window}
 A_N-1<C(N)<A_N,\qquad\lambda=C(N).
\end{equation}
The distinction is substantive. With
$\nu=\sqrt{A_N-C(N)}\in(0,1)$, the radial minimal continuum operator
has distinct nonnegative Friedrichs and Krein realizations
\cite[Proposition~4.19]{BDG11}. Consistency on tests supported away from
the origin and nonnegativity therefore do not select the Friedrichs
boundary behavior by themselves. A possible route to resolving this endpoint is to identify the asymptotic
behavior of the critical lattice state and the continuum boundary condition
it selects. By Lemma~\ref{cont:endpoint-dimensions} and \eqref{eq:known-transition}, the only dimension not covered by our endpoint results is $N=9$.
We do not assert a continuum-limit identification in
\eqref{cont:unresolved-window}.

For $\lambda>A_N$, the continuum nonexistence theorem of Baras and
Goldstein \cite{BG84} and the mesh-refinement blow-up point in the same
direction. For $C(N)<\lambda\le A_N$, Corollary~\ref{cont:mismatch}
instead shows that this sampled singular discretization fails to
approximate a well-defined continuum evolution. The lattice threshold separates uniform contractivity from mesh-refinement
blow-up. Convergence at the remaining critical endpoint needs a separate
identification.
\section{The explicit Green operator}\label{gr:green}\label{sec:green}

\subsection{The full-lattice heat kernel}

We first derive the kernel with the normalization of the Laplacian used
throughout this paper. The resulting Bessel representation agrees with the
one-dimensional heat kernel in \cite{Ciaurri}; the coordinate generators
commute, so their kernels multiply.

For $m\in\mathbb N_0$, the modified Bessel function is given by
\[
 I_m(2t)=\sum_{k=0}^{\infty}\frac{t^{2k+m}}{k!(k+m)!},
 \qquad I_{-m}(2t)=I_m(2t).
\]
Let $L_N^{\mathrm{full}}$ be the positive Laplacian on $\ell^2(\Z^N)$.
Expanding the exponentials of its commuting coordinate shifts gives
\begin{equation}\label{gr:heat}
 p_t(z):=\langle\delta_0,e^{-tL_N^{\mathrm{full}}}\delta_z\rangle
 =e^{-2Nt}\prod_{j=1}^N I_{|z_j|}(2t),\qquad t>0.
\end{equation}
Indeed, in one coordinate the coefficient of displacement $m\ge0$ in
$e^{t(S+S^{-1})}$ is the sum of $t^{2k+m}/(k!(k+m)!)$ over $k\ge0$.
The factor $e^{-2Nt}$ accounts for the diagonal $2N$.
In particular, $p_t(z)>0$ and $\sum_zp_t(z)=1$.
The Fourier representation is
\begin{equation}\label{gr:fourier}
 p_t(z)=\frac1{(2\pi)^N}\int_{[-\pi,\pi]^N}
 e^{-t\ell_N(\xi)}e^{-iz\cdot\xi}\dd\xi,
 \qquad \ell_N(\xi)=2\sum_{j=1}^N(1-\cos\xi_j).
\end{equation}
Since $1-\cos s\ge2s^2/\pi^2$ for $|s|\le\pi$, comparison with a
Gaussian integral yields
\begin{equation}\label{gr:heat-bound}
 0<p_t(z)\le p_t(0)\le B_N(1+t)^{-N/2}.
\end{equation}
For $0<t\le1$ one uses the total mass bound. Hence, for $N\ge3$,
\begin{equation}\label{gr:full-green}
 g_N(z):=\int_0^\infty p_t(z)\dd t
 =\int_0^\infty e^{-2Nt}\prod_{j=1}^N I_{|z_j|}(2t)\dd t
\end{equation}
is finite and strictly positive, with $g_N(z)\le g_N(0)$.

\subsection{The Dirichlet correction at one vertex}

We obtain the correction algebraically at positive resolvent parameter.
This is the infinite-lattice counterpart of the one-vertex Dirichlet
construction for finite graphs in \cite[Theorem~0.67]{KLW21}. The argument
below establishes the formula directly in our setting and specifies
both the boundary condition and the normalization. The corresponding probabilistic identity
follows by stopping the walk at its first visit to the origin; see
\cite{LawlerLimic}.

For $s>0$, set
\[
 g_{N,s}(z)=\int_0^\infty e^{-st}p_t(z)\dd t.
\]
This is the kernel of $(L_N^{\mathrm{full}}+s)^{-1}$.
For a finitely supported function $a$ on $\Om$, extend $a$ by zero to
$\Z^N$, and denote the extension by $\widetilde a$. Define
\[
 \alpha=\frac{\bigl((L_N^{\mathrm{full}}+s)^{-1}
                \widetilde a\bigr)(0)}{g_{N,s}(0)},
 \qquad
 v=(L_N^{\mathrm{full}}+s)^{-1}
                 (\widetilde a-\alpha\delta_0).
\]
The denominator is positive and the choice of $\alpha$ gives $v(0)=0$.
At every $x\ne0$,
$(L_N^{\mathrm{full}}+s)v(x)=a(x)$. Since $v$ vanishes at the origin,
this is exactly $(L_N+s)v|_{\Om}=a$ for our Dirichlet restriction.
Uniqueness of the $\ell^2$ resolvent solution therefore proves
\begin{equation}\label{gr:positive-green}
 (L_N+s)^{-1}(x,y)
 =g_{N,s}(x-y)-\frac{g_{N,s}(x)g_{N,s}(y)}{g_{N,s}(0)}.
\end{equation}

Let $\operatorname{Adj}_{\Om}$ denote adjacency on the punctured lattice.
The equality $L_N=2NI-\operatorname{Adj}_{\Om}$ implies
\[
 e^{-tL_N}=e^{-2Nt}\sum_{k=0}^{\infty}
                  \frac{t^k}{k!}\operatorname{Adj}_{\Om}^{k}.
\]
Its entries count paths that avoid zero and are nonnegative; they are
strictly positive for $t>0$ because $\Om$ is connected.
Deleting paths can only decrease their count, so these entries are bounded
above by $p_t(x-y)$. We may consequently integrate them in time and let
$s\downarrow0$ in \eqref{gr:positive-green}. Monotone convergence on the
resolvent integral and ordinary convergence in the explicit expression give
\begin{equation}\label{gr:killed-green}
 G_N^D(x,y):=\int_0^\infty(e^{-tL_N})(x,y)\dd t
 =g_N(x-y)-\frac{g_N(x)g_N(y)}{g_N(0)}.
\end{equation}
In particular,
\[
 0<G_N^D(x,y)\le g_N(x-y)<\infty.
\]
This argument also justifies positivity of the difference on the right of
\eqref{gr:killed-green}.

Finite Dirichlet domains provide the approximation needed to identify the weighted inverse.

\begin{lemma}\label{gr:exhaustion}\label{lem:greenlimit}
Let $F_m$ be increasing nonempty finite subsets of $\Om$ whose union is $\Om$.
Let $L_{F_m}$ be the Dirichlet restriction with diagonal $2N$, and extend
its inverse by zero when one index lies outside $F_m$. Then
\begin{equation}\label{gr:exhaustion-limit}
 L_{F_m}^{-1}(x,y)\uparrow G_N^D(x,y)
 \qquad(x,y\in\Om).
\end{equation}
\end{lemma}
\begin{proof}
If $u$ is supported in a finite set and has zero energy, all its edge
differences vanish. Connectedness of the full lattice and vanishing
outside that set imply $u=0$. Thus $L_F$ is positive definite. The symmetric nonnegative adjacency
matrix has norm equal to its largest eigenvalue, which is strictly less
than $2N$ by this positive definiteness. Writing
$P_F=\operatorname{Adj}_F/(2N)$, its norm-convergent Neumann series gives
\[
 L_F^{-1}=\frac1{2N}\sum_{k=0}^\infty P_F^k.
\]
For fixed $x,y\in\Om$ and $k$, every path of length $k$ from $x$ to $y$ in $\Om$ lies in $F_m$ for all sufficiently large $m$. Indeed, there are only finitely many such paths. The entries of $P_{F_m}^k$ therefore increase to the
corresponding punctured-lattice path counts divided by $(2N)^k$.
Monotone convergence in $m$ and $k$ gives the series obtained by
integrating the preceding heat-kernel expansion, namely $G_N^D$.
This is the concrete exhaustion principle underlying
\cite[Theorem~6.26]{KLW21}.
\end{proof}

\subsection{Identification with the inverse Hardy operator}\label{sec:inverse}

Define
\begin{equation}\label{gr:T-kernel}
 T_N(x,y)=\frac{G_N^D(x,y)}{|x||y|},\qquad x,y\in\Om.
\end{equation}
The form core identifies this kernel with the bounded inverse of $\K$
and fixes its boundary condition at infinity.

\begin{proposition}\label{gr:inverse}\label{prop:inverse}
The kernel \eqref{gr:T-kernel} defines the bounded positive selfadjoint
operator $\T=\K^{-1}$ on $\ell^2(\Om)$. Every kernel entry is strictly
positive, and
\begin{equation}\label{gr:norm-formula}
 C(N)=\frac1{\|\T\|}.
\end{equation}
\end{proposition}
\begin{proof}
Let $F_m$ be as in Lemma~\ref{gr:exhaustion}, and let
$\iota_m:\ell^2(F_m)\to\ell^2(\Om)$ be extension by zero.
The restriction of the form $\mathfrak q_N$ to this finite-dimensional
space has matrix $K_m=D_mL_{F_m}D_m$, where
$D_m=\operatorname{diag}(|x|:x\in F_m)$. Consequently,
\[
 K_m^{-1}=D_m^{-1}L_{F_m}^{-1}D_m^{-1}.
\]
We claim that $R_m:=\iota_mK_m^{-1}\iota_m^*$ converges strongly to
$\K^{-1}$. To prove this, fix $a\in\ell^2(\Om)$ and put
$z=\K^{-1}a$, $z_m=R_ma$. Their form equations are
\[
 \mathfrak q_N(z,v)=\langle a,v\rangle\quad(v\in\mathcal Q_N),
 \qquad
 \mathfrak q_N(z_m,v)=\langle a,v\rangle
            \quad(\operatorname{supp}v\subset F_m).
\]
Thus $z_m$ is the orthogonal projection of $z$ onto the functions supported
in $F_m$, in the Hilbert space $\mathcal Q_N$ with inner product
$\mathfrak q_N$. This is a Hilbert space because coercivity makes the
energy norm equivalent to its form norm. The union of these subspaces
contains $c_{00}(\Om)$ and is a form core. Orthogonal projection therefore
gives $\mathfrak q_N[z-z_m]\to0$. Coercivity then gives $z_m\to z$ in
$\ell^2$.

Applying this convergence to coordinate vectors and using
Lemma~\ref{gr:exhaustion}, we obtain
\[
 \langle\delta_x,\K^{-1}\delta_y\rangle
 =\lim_m\frac{L_{F_m}^{-1}(x,y)}{|x||y|}
 =\frac{G_N^D(x,y)}{|x||y|}.
\]
Hence the matrix entries of the bounded operator $\K^{-1}$ are exactly
\eqref{gr:T-kernel}. Each row belongs to $\ell^2$, so the series
\[
 (\T f)(x)=\sum_{y\in\Om}T_N(x,y)f(y)
\]
converges absolutely for every $f\in\ell^2$, by Cauchy--Schwarz, and
represents the operator action. Selfadjointness and positivity follow
from $\T=\K^{-1}$, and strict positivity of the entries follows from
\eqref{gr:killed-green}. Finally, \eqref{eq:bottom} and the spectral
theorem yield $\|\T\|=1/\inf\sigma(\K)=1/C(N)$.
\end{proof}

Finite compressions of the Green operator then approximate its norm monotonically.

\begin{corollary}\label{gr:compressions}
Let $H_m=(T_N(x,y))_{x,y\in F_m}$. Then
\[
 \lambda_{\max}(H_m)\uparrow\|\T\|,
 \qquad
 \lambda_{\max}(H_m)^{-1}\downarrow C(N).
\]
\end{corollary}
\begin{proof}
The largest eigenvalue is the supremum of the Rayleigh quotient of $\T$
over vectors supported in $F_m$. These spaces increase and their union
is dense in $\ell^2$. Boundedness of $\T$ gives convergence to its largest
spectral value, which is its norm because $\T$ is positive. This value
and all the displayed eigenvalues are positive, so reciprocals converge
in the stated direction.
\end{proof}

The matrix $H_m$ uses the Green kernel of the entire punctured lattice.
The matrix $K_m^{-1}$ instead uses the kernel killed on leaving $F_m$.
The distinction is relevant when using finite domains numerically:
compression of the explicit Green kernel retains excursions outside
$F_m$.

\section{Two-sided bounds from the Green operator}\label{gr:bounds}\label{sec:bounds}

We now apply the standard Rayleigh principle and weighted Schur test to the explicit Green kernel. The resulting two envelopes recover $C(N)$ exactly; the variational principles themselves are not new.

\begin{proposition}\label{gr:envelopes}
For $f\in\ell^2(\Om)$ with $f(x)>0$ at every vertex, put
\[
 q_N(f)=\frac{\langle f,\T f\rangle}{\|f\|_2^2},
 \qquad
 s_N(f)=\sup_{x\in\Om}\frac{(\T f)(x)}{f(x)}.
\]
Then $q_N(f)>0$, $0<s_N(f)\le\infty$, and
\begin{equation}\label{gr:two-sided}
 \frac1{s_N(f)}\le C(N)\le\frac1{q_N(f)},
 \qquad \frac1\infty:=0.
\end{equation}
Both envelopes are exact:
\begin{equation}\label{gr:exact-envelopes}
 C(N)=\sup_{\substack{f\in\ell^2(\Om)\\f>0}}\frac1{s_N(f)}
     =\inf_{\substack{f\in\ell^2(\Om)\\f>0}}\frac1{q_N(f)}.
\end{equation}
\end{proposition}
\begin{proof}
Write $T=\T$ and $\rho=\|T\|$. Since $T=\K^{-1}$, the spectral theorem
gives $\langle f,Tf\rangle>0$ for every nonzero $f$.
The Rayleigh principle gives $q_N(f)\le\rho$.

If $s=s_N(f)<\infty$ and $v$ is finitely supported, the symmetry and
nonnegative entries of $T$ imply
\begin{align*}
 |\langle v,Tv\rangle|
 &\le\sum_{x,y\in\operatorname{supp}v}T_N(x,y)|v(x)||v(y)|\\
 &\le\frac12\sum_{x,y\in\operatorname{supp}v}T_N(x,y)
 \left(|v(x)|^2\frac{f(y)}{f(x)}
       +|v(y)|^2\frac{f(x)}{f(y)}\right)\\
 &\le\sum_x|v(x)|^2\frac{(Tf)(x)}{f(x)}
 \le s\|v\|_2^2.
\end{align*}
All extensions of sums in the third line involve nonnegative terms.
Density and selfadjointness give $\rho\le s$.
Taking reciprocals proves \eqref{gr:two-sided}; the case $s=\infty$
follows with the indicated convention.

For exactness of the Rayleigh envelope, positivity of the kernel gives
$\langle v,Tv\rangle\le\langle |v|,T|v|\rangle$.
The supremum of the Rayleigh quotient can therefore be taken over
nonnegative vectors. Fix $g\in\ell^2$ with $g>0$ everywhere.
Approximating any nonnegative vector $v$ by $v+\varepsilon g$ shows that
strictly positive vectors have the same supremum. Hence
$\sup_{f>0}q_N(f)=\rho$.

For the Schur envelope, let $a>\rho$ and define
\[
 f_a=(aI-T)^{-1}g
 =\frac1a\sum_{j=0}^{\infty}(T/a)^jg.
\]
The series converges in $\ell^2$, all its terms are nonnegative, and its
first term is strictly positive. Thus $f_a>0$, and
$Tf_a=af_a-g\le af_a$ gives $s_N(f_a)\le a$.
Consequently $\sup_{f>0}1/s_N(f)\ge1/a$.
Letting $a\downarrow\rho$ proves the remaining equality.
\end{proof}

These bounds can be evaluated through the physical profile $u$. The final local bound is the usual positive-supersolution consequence of the ground-state identity, written for the prescribed weight $|x|^{-2}$.

\begin{proposition}\label{gr:profiles}\label{eff:local-profile}
Let $u>0$ on $\Om$, set $u(0)=0$, and assume
$\E[u]+\M[u]<\infty$. Define
\[
 V_u(x)=\sum_{y\in\Om}G_N^D(x,y)\frac{u(y)}{|y|^2},
 \qquad
 Q_G[u]=\sum_{x,y\in\Om}G_N^D(x,y)
       \frac{u(x)u(y)}{|x|^2|y|^2}.
\]
These expressions are finite, and
\begin{equation}\label{gr:profile-bounds}
 \left[\sup_x\frac{V_u(x)}{u(x)}\right]^{-1}
 \le C(N)\le\frac{\M[u]}{Q_G[u]}
 \le\frac{\E[u]}{\M[u]}.
\end{equation}
Moreover, without any integrability assumption on a positive $u$, put
\begin{equation}\label{eff:Hu}
 H_u(x)=|x|^2\frac{L_Nu(x)}{u(x)}.
\end{equation}
If $H_u(x)\ge\gamma$ at every vertex, then $C(N)\ge\gamma$.
The latter assertion also holds on every Dirichlet subdomain.
\end{proposition}
\begin{proof}
Put $f=D^{-1}u$. Lemma~\ref{lem:cutoff} gives $f\in\mathcal Q_N$,
$\|f\|_2^2=\M[u]$, and $\mathfrak q_N[f]=\E[u]$.
The kernel representation yields
\[
 (\T f)(x)=\frac{V_u(x)}{|x|},\qquad
 \langle f,\T f\rangle=Q_G[u].
\]
The pointwise series is absolutely convergent; Tonelli's theorem and
boundedness of $\T$ justify the double sum. Proposition~\ref{gr:envelopes}
gives the first two inequalities in \eqref{gr:profile-bounds}.
For the third, the spectral theorem and Cauchy--Schwarz give
\[
 \M[u]^2=|\langle\K^{1/2}f,\K^{-1/2}f\rangle|^2
 \le\E[u]Q_G[u].
\]
Finally, the ground-state identity, Lemma~\ref{lem:gs}, applied to an
arbitrary finitely supported test function $v$, gives
\[
 \E[v]\ge\sum_{x\in\Om}\frac{L_Nu(x)}{u(x)}|v(x)|^2
 \ge\gamma\M[v].
\]
This calculation is finite for every such $v$, so it requires no
integrability assumption on $u$. The same edge calculation, with zero
values on the complement, proves the Dirichlet-subdomain assertion.
\end{proof}

\section{The growth rate from the shifted Green family}\label{gr:shifted}

The same Green construction determines the supercritical growth rate.
Write $W=V^{1/2}$, so that $(Wf)(x)=f(x)/|x|$, and for $s>0$ set
\begin{equation}\label{gr:Ts}
 \T(s)=W(L_N+s)^{-1}W.
\end{equation}
Formula \eqref{gr:positive-green} makes its kernel explicit:
\begin{equation}\label{gr:Ts-kernel}
 T_N(s;x,y)=\frac1{|x||y|}
 \left(g_{N,s}(x-y)
       -\frac{g_{N,s}(x)g_{N,s}(y)}{g_{N,s}(0)}\right).
\end{equation}

To locate the growth rate, we first control the norm of this family from zero energy to infinity.

\begin{proposition}\label{gr:Ts-limit}
For every $s>0$, $\T(s)$ is nonzero, compact, positive and selfadjoint.
As $s\downarrow0$,
\begin{equation}\label{gr:Ts-strong}
 \T(s)\longrightarrow\T\quad\hbox{strongly},
 \qquad \|\T(s)\|\longrightarrow\|\T\|=1/C(N).
\end{equation}
The function $s\mapsto\|\T(s)\|$ is continuous and strictly decreasing
on $(0,\infty)$, and tends to zero as $s\to\infty$.
\end{proposition}
\begin{proof}
The diagonal multiplication operator $W$ is compact, since $|x|^{-1}$
tends to zero at infinity: truncating $W$ to a finite ball converges in
operator norm. Thus \eqref{gr:Ts} is compact. Positivity and selfadjointness
follow from those of $(L_N+s)^{-1}$; strict positivity of its diagonal
shows that it is nonzero.

We first identify its endpoint without assuming a bounded inverse for
$L_N$. The variational identity for the positive invertible operator
$L_N+s$ is
\begin{equation}\label{gr:resolvent-variational}
 \langle f,\T(s)f\rangle
 =\sup_{u\in\ell^2(\Om)}
 \left(2\operatorname{Re}\langle Wf,u\rangle
       -\E[u]-s\|u\|_2^2\right).
\end{equation}
It follows by completing the square in $(L_N+s)^{1/2}u$.
For every $u\in\ell^2$, both its energy and its weighted mass are finite,
because $L_N$ and $W$ are bounded. Lemma~\ref{lem:cutoff} therefore shows
that $z=Wu\in\mathcal Q_N$, with $\mathfrak q_N[z]=\E[u]$.
Dropping the last negative term in \eqref{gr:resolvent-variational} and
enlarging the admissible set to $\mathcal Q_N$ yields
\[
 \langle f,\T(s)f\rangle
 \le\sup_{z\in\mathcal Q_N}
       \left(2\operatorname{Re}\langle f,z\rangle
                          -\mathfrak q_N[z]\right)
 =\langle f,\K^{-1}f\rangle.
\]
The last equality follows by completing the square in $\K^{1/2}z$.
Conversely, each finitely supported $z$ can be written as $z=Wu$ with
$u=Dz$ finitely supported. For this fixed $u$, the expression in
\eqref{gr:resolvent-variational} tends to
$2\operatorname{Re}\langle f,z\rangle-\mathfrak q_N[z]$ as $s\downarrow0$.
The form-core property thus gives
\[
 \langle f,\T(s)f\rangle\uparrow\langle f,\T f\rangle
 \qquad(f\in\ell^2).
\]
Here monotonicity also follows directly from the resolvent identity.
Since $0\le\T-\T(s)\le\T$, we have
\[
 \|(\T-\T(s))f\|^2
 \le\|\T\|\langle f,(\T-\T(s))f\rangle\longrightarrow0.
\]
This proves strong convergence. The upper bound
$\|\T(s)\|\le\|\T\|$ and the convergence of every unit-vector Rayleigh
quotient imply convergence of the norms.

For $s,t>0$, the resolvent identity gives
\[
 \|\T(s)-\T(t)\|\le\frac{|s-t|}{st},
\]
so the family is norm-continuous on $(0,\infty)$.
To prove strict decrease, take $0<s<t$ and a unit eigenvector $f$ of
$\T(t)$ for its largest eigenvalue $\|\T(t)\|$; such an eigenvector
exists by compactness and positivity. Since $W$ is injective, $Wf\ne0$.
The spectral theorem for $L_N$ gives
\begin{align*}
 \langle f,(\T(s)-\T(t))f\rangle
 &=(t-s)\langle Wf,(L_N+s)^{-1}(L_N+t)^{-1}Wf\rangle\;>0.
\end{align*}
Consequently $\|\T(s)\|\ge\langle f,\T(s)f\rangle>\|\T(t)\|$.
Finally, $\|\T(s)\|\le\|W\|^2/s=1/s$ proves the limit at infinity.
\end{proof}

The Birman--Schwinger factorization converts this monotonicity into a unique scalar equation for the growth rate.

\begin{proposition}\label{gr:BS}
Let $\lambda>C(N)$ and let
$\mu_N(\lambda)=-\inf\sigma(L_N-\lambda V)>0$.
Then $\mu_N(\lambda)$ is the unique positive solution of
\begin{equation}\label{gr:mu-equation}
 \lambda\|\T(s)\|=1.
\end{equation}
More generally, for every $s>0$ and $\lambda>0$,
\begin{equation}\label{gr:BS-positive}
 L_N-\lambda V+s\ge0
 \quad\Longleftrightarrow\quad\lambda\|\T(s)\|\le1,
\end{equation}
and $-s$ is an eigenvalue of $L_N-\lambda V$ if and only if
$1/\lambda$ is an eigenvalue of $\T(s)$.
\end{proposition}
\begin{proof}
Let $A_s=W(L_N+s)^{-1/2}$. The elementary identity
$\|A_s^*A_s\|=\|A_sA_s^*\|=\|A_s\|^2$ and
$A_sA_s^*=\T(s)$ show that
\[
 L_N-\lambda V+s
 =(L_N+s)^{1/2}(I-\lambda A_s^*A_s)(L_N+s)^{1/2}
\]
is nonnegative exactly when $\lambda\|\T(s)\|\le1$.
By Proposition~\ref{gr:Ts-limit}, the norm decreases continuously and
strictly from $1/C(N)$ to zero. Since $\lambda>C(N)$, equation
\eqref{gr:mu-equation} has a unique root $s_*$. Equivalence
\eqref{gr:BS-positive} says that $L_N-\lambda V+s$ is nonnegative
exactly for $s\ge s_*$. Its spectral bottom is therefore $-s_*$,
which proves $s_*=\mu_N(\lambda)$.

For the eigenvalue assertion, suppose
$(L_N-\lambda V)\phi=-s\phi$, $\phi\ne0$, and set $v=W\phi$.
Then $v\ne0$ and
\[
 v=\lambda W(L_N+s)^{-1}Wv=\lambda\T(s)v.
\]
Conversely, if $\T(s)v=\lambda^{-1}v$, set
$\phi=(L_N+s)^{-1}Wv$.
The equality $W\phi=\lambda^{-1}v$ shows that $\phi\ne0$ and
\[
 (L_N+s)\phi=Wv=\lambda V\phi.
\]
Thus $(L_N-\lambda V)\phi=-s\phi$.
\end{proof}

The zero-energy member of this family gives the critical coupling
$C(N)$, while its positive-parameter members determine the rate of the
supercritical heat evolution.
\section{Symmetry reduction and an explicit exterior gap}\label{eff:gap-section}

The orbit reduction and the first-sphere Schur complement used in this section were developed for the three-term expansion in \cite{LizamaAsymptotics}. We retain that framework and derive the explicit exterior bound needed for the uniform remainder estimate below.

Write $s(x)=|x|^2$, and let $\mathfrak G_N$ be the finite group generated
by coordinate permutations and independent changes of coordinate signs.
The operator $\K=DL_ND$, where $D=\sqrt{s}$, preserves the closed subspace
$\mathcal H_N^{\mathrm{inv}}\subset\ell^2(\Om)$ of
$\mathfrak G_N$-invariant functions. In this section, a Dirichlet restriction
always retains the diagonal $2N$ of $L_N$.

The symmetry allows us to retain the sharp constant in the invariant subspace.

\begin{lemma}\label{eff:symmetry}
The infimum of the spectrum of the restriction of $\K$ to
$\mathcal H_N^{\mathrm{inv}}$ is $C(N)$.
\end{lemma}
\begin{proof}
For a finitely supported physical function $u$, extended by zero to the
origin, define
\[
 u_{\mathrm{av}}(x)=
 \left(\frac1{|\mathfrak G_N|}
       \sum_{g\in\mathfrak G_N}|u(gx)|^2\right)^{1/2}.
\]
This function is invariant and finitely supported. Since the group preserves
$s(x)$, a change of variables gives
$\M[u_{\mathrm{av}}]=\M[u]$. For each edge $x\sim y$, apply the reverse
triangle inequality to the two vectors
$|\mathfrak G_N|^{-1/2}(u(gx))_g$ and
$|\mathfrak G_N|^{-1/2}(u(gy))_g$. It gives
\[
 |u_{\mathrm{av}}(x)-u_{\mathrm{av}}(y)|^2
 \le \frac1{|\mathfrak G_N|}
       \sum_{g\in\mathfrak G_N}|u(gx)-u(gy)|^2.
\]
Summing over unoriented lattice edges, which the group permutes, proves
$\E[u_{\mathrm{av}}]\le\E[u]$. Thus the infimum in
\eqref{eq:constant} is unchanged by restricting to invariant functions.
Multiplication by $D$ preserves invariance, so the form-core construction
of $\K$ and \eqref{eq:bottom} prove the assertion.
\end{proof}

Let
\[
 Q_N=\{x\in\Z^N:s(x)\ge2\},
 \qquad
 L_{N,Q}v(x)=2Nv(x)-\sum_{\substack{y\sim x\\y\in Q_N}}v(y).
\]
Let $\mathcal K_{N,Q}$ be the Friedrichs realization of
$\sqrt{s}\,L_{N,Q}\sqrt{s}$ in $\ell^2(Q_N)$.

An explicit lower bound on the exterior will separate the first sphere from the remaining spectral problem.

\begin{lemma}\label{eff:gap}
For every integer $N\ge44$,
\begin{equation}\label{eff:gap-bound}
 \mathcal K_{N,Q}\ge\frac{20N}{9}I
\end{equation}
as a quadratic-form inequality. In particular, the same lower bound holds
after restriction to invariant functions.
\end{lemma}
\begin{proof}
Set $h(x)=s(x)^{-2}$ on $Q_N$ and extend $h$ by zero to its complement.
The edge calculation in Lemma~\ref{lem:gs}, now on $Q_N$, gives for every
finitely supported $v$ on $Q_N$,
\begin{align}\label{eff:gap-gs}
 \langle\mathcal K_{N,Q}v,v\rangle
 &=\sum_{\substack{\{x,y\}\subset Q_N\\x\sim y}}
 h(x)h(y)
 \left|\frac{\sqrt{s(x)}v(x)}{h(x)}
       -\frac{\sqrt{s(y)}v(y)}{h(y)}\right|^2
 +\sum_{x\in Q_N}\mathcal H_h(x)|v(x)|^2,\\
 \mathcal H_h(x)&=s(x)\frac{L_{N,Q}h(x)}{h(x)}.\notag
\end{align}
No summability assumption on $h$ is needed, since $v$ has finite support.
It remains to prove $\mathcal H_h\ge20N/9$ pointwise.

Suppose first that $s=s(x)\ge6$. Every nearest neighbor of $x$ belongs
to $Q_N$, since its squared norm is at least
$(\sqrt{s}-1)^2>2$. Put $a=s+1$ and $t_i=x_i^2$, so that
$\sum_{i=1}^Nt_i=s$. For $0\le t\le s$,
\begin{align*}
 f_s(t)&=(a+2\sqrt t)^{-2}+(a-2\sqrt t)^{-2}\\
 &=\frac{2(a^2+4t)}{(a^2-4t)^2}
 =\frac{2}{a^2}\sum_{k=0}^\infty(2k+1)
       \left(\frac{4t}{a^2}\right)^k.
\end{align*}
The series converges uniformly on $[0,s]$, because $4s/(s+1)^2<1$.
For every integer $k\ge1$, $\sum_i t_i^k\le s^k$. Its nonnegative
coefficients therefore imply
\[
 \sum_{i=1}^N f_s(t_i)\le (N-1)f_s(0)+f_s(s).
\]
Since $h(x)^{-1}=s^2$, this gives
\begin{equation}\label{eff:gap-radial}
 \mathcal H_h(x)
 =2Ns-s^3\sum_{i=1}^N f_s(t_i)
 \ge\alpha(s)N-B(s),
\end{equation}
where
\[
 \alpha(s)=\frac{2s(2s+1)}{(s+1)^2},
 \qquad
 B(s)=\frac{8s^4(3s^2+2s+3)}{(s-1)^4(s+1)^2}.
\]
Differentiating these rational functions on $s>1$ yields
\[
 \alpha'(s)=\frac{2(3s+1)}{(s+1)^3}>0,
 \qquad
 B'(s)=-\frac{32s^3(s+3)(2s^2+s+1)}{(s-1)^5(s+1)^3}<0.
\]
Consequently, for $s\ge6$,
\[
 \mathcal H_h(x)\ge\frac{156}{49}N-\frac{1275264}{30625}.
\]
The latter quantity is at least $20N/9$ whenever
\begin{equation}\label{eff:gap-large-threshold}
 N\ge
 \frac{B(6)}{\alpha(6)-20/9}
 =\frac{1434672}{33125}<44.
\end{equation}

It remains to check $s=2,3,4,5$. Up to signed permutations, the nonzero
coordinates form one of the six patterns in the following table. The
last column is a sufficient lower bound on $N$ for the middle column to
be at least $20N/9$.
\begin{center}
\renewcommand{\arraystretch}{1.35}
\begin{tabular}{@{}ccc@{}}
\toprule
Pattern & $\mathcal H_h(x)$ & Sufficient bound on $N$\\
\midrule
$(1,1)$ & $20N/9+656/225$ & $0$\\
$(1,1,1)$ & $21N/8-99/8$ & $891/29$\\
$(2)$ & $72N/25+8768/2025$ & $0$\\
$(1,1,1,1)$ & $72N/25-145408/11025$ & $36352/1813$\\
$(2,1)$ & $55N/18-16345/576$ & $3269/96$\\
$(1,1,1,1,1)$ & $55N/18-8125/576$ & $1625/96$\\
\bottomrule
\end{tabular}
\end{center}
For completeness, every entry follows from the same neighbor formula.
If the pattern is $p=(p_1,\ldots,p_\ell)$ and
$s=\sum_{j=1}^{\ell}p_j^2$, then
\begin{align}\label{eff:gap-small-formula}
 \mathcal H_h(p)=2Ns-s^3\left[
 \frac{2(N-\ell)}{(s+1)^2}
 +\sum_{j=1}^{\ell}\left(
 \frac{\mathbf1_{\{s+1-2p_j\ge2\}}}{(s+1-2p_j)^2}
 +\frac1{(s+1+2p_j)^2}\right)\right].
\end{align}
For example, for $(1,1)$ the neighbor sum is
$2/25+2(N-2)/9$, and for $(2)$ it is
$1/81+2(N-1)/25$. The neighbors with squared norm one have been
removed in both cases, while the diagonal $2N$ is retained. Substitution
gives the first and third rows. The other four substitutions respectively
use the neighbor sums
\[
 \begin{gathered}
 \frac34+\frac3{36}+\frac{2(N-3)}{16},\qquad
 \frac49+\frac4{49}+\frac{2(N-4)}{25},\\
 \frac14+\frac1{100}+\frac1{16}+\frac1{64}
       +\frac{2(N-2)}{36},\qquad
 \frac5{16}+\frac5{64}+\frac{2(N-5)}{36}.
 \end{gathered}
\]
All thresholds in the table are less than $44$. Combining these
calculations with \eqref{eff:gap-large-threshold} proves the pointwise
bound. Equation~\eqref{eff:gap-gs} proves \eqref{eff:gap-bound} on the
form core, and closure extends it to the full form domain.
\end{proof}

This separation reduces the spectral problem to the first sphere and an invertible exterior block.

\begin{lemma}\label{eff:first-shell}
Let $N\ge44$. In the invariant subspace, separate the normalized indicator
$\varphi_1=(2N)^{-1/2}\mathbf1_{\{s=1\}}$. Then
\begin{equation}\label{eff:first-shell-block}
 \K\big|_{\mathcal H_N^{\mathrm{inv}}}
 =\begin{pmatrix}2N&-b_N^*\\-b_N&H_N\end{pmatrix},
 \qquad H_N\ge\frac{20N}{9}I.
\end{equation}
In the orthonormal orbit basis of the exterior, the only nonzero
coordinates of $b_N$ are
\begin{equation}\label{eff:first-shell-source}
 b_N(1,1)=\sqrt{8(N-1)},\qquad b_N(2)=2,
 \qquad \|b_N\|^2=8N-4.
\end{equation}
The number $C(N)$ is a simple isolated eigenvalue of this invariant
restriction, lies below $2N$, and is its only spectral point below $20N/9$.
\end{lemma}
\begin{proof}
The set $\{s=1\}$ is a single orbit, so its invariant subspace is
one-dimensional. The exterior diagonal block is precisely the restriction
of $\mathcal K_{N,Q}$ to invariant functions, giving its lower bound by
Lemma~\ref{eff:gap}.

The first sphere contains $2N$ vertices. The orbits $(1,1)$ and $(2)$
contain $2N(N-1)$ and $2N$ vertices, respectively. There are
$4N(N-1)$ edges from the first sphere to $(1,1)$ and $2N$ edges to
$(2)$. The corresponding edge weights of $\sqrt{s}\,\Adj\sqrt{s}$
are $\sqrt2$ and $2$. Dividing each weighted edge count by the square
root of the product of the two orbit sizes gives
\eqref{eff:first-shell-source}. These are the only exterior orbits
adjacent to the first sphere.

Let $\varphi_{11}$ be the normalized indicator of the orbit $(1,1)$.
For $t>0$,
\[
 \frac{\langle\K(\varphi_1+t\varphi_{11}),
                     \varphi_1+t\varphi_{11}\rangle}{1+t^2}
 =\frac{2N-2t\sqrt{8(N-1)}+4Nt^2}{1+t^2}<2N
\]
when $t$ is sufficiently small. Lemma~\ref{eff:symmetry} therefore gives
$C(N)<2N$.
Finally, \eqref{eff:first-shell-block} is a bounded finite-rank
perturbation of $2N\oplus H_N$. Its essential spectrum is contained
in $[20N/9,\infty)$, so its bottom below $2N$ is an isolated
eigenvalue. There can be at most one eigenvalue, counting multiplicity,
below $20N/9$: a two-dimensional spectral subspace below this threshold
would contain a nonzero vector orthogonal to $\varphi_1$, whose Rayleigh
quotient is at least $20N/9$ by the exterior bound. This proves simplicity
and the claimed separation.
\end{proof}

\section{An effective rational approximation in large dimension}\label{eff:error-section}

The Green representation gives exact two-sided variational formulas for every
integer dimension under consideration. We now obtain an explicit approximation
whose error tends to zero with a specified rate and threshold. The construction
uses the first sphere, a finite collection of return paths, and an analytic
estimate for their complement.

We prove Theorem~\ref{eff:main} for the rational function
$R(N)$ defined in \eqref{eff:R-intro}.

The proof below uses the exterior estimate of Lemma~\ref{eff:gap} and the
asymptotic expansion established in Proposition~\ref{eff:asym-expansion}.
Every estimate in the remainder argument is uniform for $N\ge50$.

\subsection{Choice of the rational function}

The first three terms below recover \cite[Theorem~1.1]{LizamaAsymptotics}. The further coefficients are established in Appendix~\ref{eff:asym-appendix}:
\[
 C(N)=2N-4-\sum_{j=1}^4\frac{d_j}{N^j}+O(N^{-5}),
\]
where
\[
 d_1=\frac{20}{3},\quad d_2=\frac{253}{9},\quad
 d_3=\frac{28409}{180},\quad d_4=\frac{477372781}{453600}.
\]
We seek the finite continued-fraction expression
\[
 R(N)=2N-4-\frac{d_1}{N-a-b/(N-c)}.
\]
Expanding its last term at infinity gives the coefficients
$d_1$, $d_1a$, $d_1(a^2+b)$ and $d_1(a^3+2ab+bc)$.
Consequently the successive choices
\[
 a=\frac{d_2}{d_1},\qquad
 b=\frac{d_3}{d_1}-a^2,\qquad
 c=\frac{d_4/d_1-a^3-2ab}{b}
\]
produce exactly the rational parameters in \eqref{eff:R-intro}.
Here $b=10609/1800>0$, so all three choices are defined.
Matching these coefficients motivates $R(N)$; the remainder estimates
below establish its uniform accuracy for the infinite-lattice constant.

\subsection{A stable scalar equation}

Use the invariant subspace and first-sphere decomposition from
Lemma~\ref{eff:first-shell}. In its orthonormal orbit basis, write
\begin{equation}\label{eff:block-recall}
 \K=\begin{pmatrix}2N&-b_N^*\\-b_N&H_N\end{pmatrix},
 \qquad \|b_N\|^2=8N-4.
\end{equation}
Here $H_N$ is the exterior compression and $b_N$ has the two nonzero entries
$\sqrt{8(N-1)}$ on the orbit $(1,1)$ and $2$ on the orbit $(2)$.
For $E\le2N$, define
\begin{equation}\label{eff:schur-functions}
 m_N(E)=\langle b_N,(H_N-E)^{-1}b_N\rangle,
 \qquad F_N(E)=2N-E-m_N(E).
\end{equation}
By Lemma~\ref{eff:gap}, $H_N-E\ge(2N/9)I$ for $N\ge44$, so these
expressions are well defined throughout this interval.

A lower bound on $|F_N^{\prime}|$ will turn the residual of $F_N$ into
an error bound for $C(N)$.

\begin{lemma}\label{eff:scalar-stability}
For $N\ge44$, $C(N)$ is the unique zero of $F_N$ below $2N$. Moreover,
\begin{equation}\label{eff:F-derivative}
 F_N'(E)=-1-\|(H_N-E)^{-1}b_N\|^2,
 \qquad
 1\le |F_N'(E)|\le1+\frac{162}{N}
 \quad(E\le2N).
\end{equation}
In particular, for every real $E\le2N$,
\begin{equation}\label{eff:scalar-error}
 |C(N)-E|\le |F_N(E)|.
\end{equation}
\end{lemma}

\begin{proof}
Completing the square in the exterior variable gives, for a scalar $t$ and a
vector $v$ in the exterior form domain,
\begin{align*}
 \langle (t,v),(\K-E)(t,v)\rangle
 &=F_N(E)|t|^2\\
 &\quad+
 \bigl\langle v-t(H_N-E)^{-1}b_N,
 (H_N-E)\bigl(v-t(H_N-E)^{-1}b_N\bigr)\bigr\rangle.
\end{align*}
Thus $\K-E\ge0$ if and only if $F_N(E)\ge0$. At a zero of $F_N$,
the vector $(1,(H_N-E)^{-1}b_N)$ is an eigenvector. The spectral bottom
on the invariant subspace equals $C(N)$ by Lemma~\ref{eff:symmetry}.
Also $F_N(2N)=-m_N(2N)<0$, whereas $F_N(E)\to+\infty$ as
$E\to-\infty$. Differentiation of the bounded resolvent proves the first
identity in \eqref{eff:F-derivative}, hence strict monotonicity and uniqueness
of the zero. The remaining estimate follows from
\[
 \|(H_N-E)^{-1}b_N\|^2
 \le\left(\frac9{2N}\right)^2(8N-4)
 =\frac{162}{N}-\frac{81}{N^2}\le\frac{162}{N}.
\]
Finally, integrate $-F_N'\ge1$ between $E$ and the zero $C(N)$.
\end{proof}

We record a bound for the rational expression on a complex disk. This will
also show that its real values lie in the interval needed above. Put
\begin{equation}\label{eff:delta-rho}
 \delta(z)=4+\frac{20}{3}
  \frac{z(1-cz)}{1-(a+c)z+(ac-b)z^2},
 \qquad \rho(z)=1-\frac z2\delta(z).
\end{equation}
At $z=1/N$, these functions satisfy
$\delta(1/N)=2N-R(N)$ and $\rho(1/N)=R(N)/(2N)$.
Since $a<5$, $c<6$, and $0<ac-b<30$, for $|z|\le1/25$ we have
\begin{align}
 |1-(a+c)z+(ac-b)z^2|
 &\ge1-\frac{11}{25}-\frac{30}{625}=\frac{64}{125},
 \label{eff:denominator-bound}\\
 |\delta(z)|&<4+\frac{31}{48}<5,
 \qquad |\rho(z)-1|<\frac1{10}.\label{eff:disk-bounds}
\end{align}
For real $N\ge50$, the correction to $4$ in $\delta(1/N)$ is positive.
Consequently,
\begin{equation}\label{eff:R-range}
 0<2N-5<R(N)<2N-4<2N.
\end{equation}

\subsection{An exact finite identity for the resolvent}

Fix $N\ge50$ and $E=R(N)$. On the exterior of the first sphere, let
$s(x)=|x|^2$ and write
\begin{equation}\label{eff:DBf}
 \begin{gathered}
 H_N-E=D_E-\mathcal A_{N,Q},\qquad
 D_E=\operatorname{diag}(2Ns-E),\\
 B=D_E^{-1/2}\mathcal A_{N,Q}D_E^{-1/2},\qquad f=D_E^{-1/2}b_N.
 \end{gathered}
\end{equation}
The identity in the first line is initially an identity on finite-support
vectors. Here $\mathcal A_{N,Q}$ is the weighted adjacency with nonnegative
matrix entries: before restriction to invariant functions, its entry on a
neighboring pair $x,y$ is $\sqrt{s(x)s(y)}$. No positive-semidefinite
operator property is asserted for this adjacency.
The diagonal $D_E$ is strictly positive. To define $B$ on the full exterior
space without an unbounded-operator domain ambiguity, put
\[
 a_E(x)=\sqrt{\frac{s(x)}{2Ns(x)-E}},\qquad
 B=M_{a_E}\Adj_QM_{a_E},
\]
where $M_{a_E}$ denotes multiplication and $\Adj_Q$ the unweighted exterior
adjacency. Since $s\ge2$ and $E\le2N$, one has $a_E^2\le1/N$;
therefore $B$ is bounded and selfadjoint, with $\|B\|\le2$.
It preserves the invariant subspace. Finite powers of $B$ applied to $f$
have finite support, so every unbounded-operator identity used below
is evaluated on vectors in the operator domain. Define
\[
 P_4(B)=I+B+B^2+B^3+B^4,\qquad y=D_E^{-1/2}P_4(B)f.
\]
The identity $(I-B)P_4(B)=I-B^5$ gives the residual
\begin{equation}\label{eff:finite-residual}
 r=b_N-(H_N-E)y=D_E^{1/2}B^5f.
\end{equation}

This residual gives an exact finite expansion with a nonnegative remainder.

\begin{lemma}\label{eff:finite-identity}
With the notation above,
\begin{equation}\label{eff:finite-split}
 m_N(E)=S_4(N)+\mathcal R_N,
 \qquad S_4(N)=\sum_{j=0}^4\langle f,B^{2j}f\rangle,
\end{equation}
where
\begin{equation}\label{eff:positive-tail}
 \mathcal R_N=
 \langle D_E^{1/2}B^5f,(H_N-E)^{-1}D_E^{1/2}B^5f\rangle\ge0.
\end{equation}
\end{lemma}

\begin{proof}
For the positive operator $J=H_N-E$ and $r=b_N-Jy$, direct expansion of
$\langle r,J^{-1}r\rangle$ yields
\[
 m_N(E)=2\langle b_N,y\rangle-\langle y,Jy\rangle
          +\langle r,J^{-1}r\rangle.
\]
The first two terms equal
\[
 \left\langle f,
 \bigl(2P_4(B)-P_4(B)(I-B)P_4(B)\bigr)f\right\rangle.
\]
A finite polynomial multiplication gives
\[
 2P_4(B)-P_4(B)(I-B)P_4(B)=\sum_{j=0}^9B^j.
\]
Both components of $f$ lie on the even layer $|x|_1=2$, whereas each
nearest-neighbor step changes the parity of $|x|_1$. Hence
$\langle f,B^{2j+1}f\rangle=0$. Equations
\eqref{eff:finite-residual} and \eqref{eff:positive-tail} finish the proof.
\end{proof}

\subsection{A uniform estimate for the exterior remainder}

Let $\Lambda_k=\{x\in\Om:|x|_1=k\}$ for $k\ge2$. The operator $B$ only
couples consecutive layers. We first bound these blocks on the full exterior
space; the estimates then hold on its invariant subspace as well.
For the unweighted adjacency block $J_{k,k+1}$, every row has at most $2N$
nonzero entries and every column has at most $k+1$. The row-column Schur
estimate therefore gives
\[
 \|J_{k,k+1}\|\le\sqrt{2N(k+1)}.
\]
For $x\in\Lambda_k$ and $y\in\Lambda_{k+1}$, the integer coordinates imply
$s(x)\ge k$ and $s(y)\ge k+1$. Since $E\le2N$,
\[
 \sqrt{\frac{s(x)}{2Ns(x)-E}}
 \sqrt{\frac{s(y)}{2Ns(y)-E}}
 \le\frac1{2N}\sqrt{\frac{k+1}{k-1}}.
\]
Multiplication by the diagonal factors in \eqref{eff:DBf} proves
\begin{equation}\label{eff:layer-bound}
 \|B_{k,k+1}\|\le\frac{k+1}{\sqrt{2N(k-1)}}.
\end{equation}
Also,
\begin{equation}\label{eff:source-real}
 \|f\|^2=\frac{8(N-1)}{4N-E}+\frac4{8N-E}
 \le4-\frac4N+\frac{2}{3N}<4.
\end{equation}

The exterior gap and path counts now bound the remainder uniformly in dimension.

\begin{lemma}\label{eff:tail-bound}
For every integer $N\ge50$,
\begin{equation}\label{eff:tail-constant}
 0\le\mathcal R_N<\frac{2616003280989}{12500N^5}.
\end{equation}
\end{lemma}

\begin{proof}
Expand $\|B^5f\|^2=\langle f,B^{10}f\rangle$ by inserting the orthogonal
layer projections between consecutive factors. Every surviving term is
indexed by a walk of ten steps which starts and ends at layer $2$ and never
goes below it. Its successive layer increments are $+1$ or $-1$. Subtracting
$2$ from the layer indices identifies these walks with Dyck paths of
semilength $5$, of which there are
\[
 \frac1{6}\binom{10}{5}=42.
\]
No such walk rises above layer $7$. For $2\le k\le6$,
\[
 \frac{(k+1)^2}{2(k-1)}\le\frac{49}{10}.
\]
Apply \eqref{eff:layer-bound} to each of its ten blocks and
\eqref{eff:source-real} to the endpoints. Summing the $42$ bounds gives
\begin{equation}\label{eff:B5-bound}
 \|B^5f\|^2<4\cdot42\left(\frac{49}{10N}\right)^5.
\end{equation}
The vector $B^5f$ is supported in $|x|_1\le7$, and thus in $s(x)\le49$.
By \eqref{eff:R-range}, $E>0$, so $D_E\le98N$ on this support. Moreover,
$\|(H_N-E)^{-1}\|\le9/(2N)$. Hence
\[
 \mathcal R_N\le\frac9{2N}\|D_E^{1/2}B^5f\|^2
 \le441\|B^5f\|^2.
\]
Combining this with \eqref{eff:B5-bound} gives
\[
 \mathcal R_N
 <\frac{36\cdot49\cdot42(49/10)^5}{N^5}
 =\frac{2616003280989}{12500N^5}.
\]
\end{proof}

\subsection{Holomorphic control of the finite moments}

Only the finite part $S_4$ will be continued in the parameter $z=1/N$.
Let $\mathcal P$ be the fixed set of integer partitions $p$ satisfying
$2\le |p|_1\le6$. Write $\ell(p)$ for its number of parts,
$m_j(p)$ for the multiplicity of $j$, and $s(p)=\sum_i p_i^2$.
A closed walk of at most eight steps starting on layer $2$ cannot leave
these layers. Thus the compression of $B$ to $\mathcal P$ computes every
moment in $S_4$ exactly.

We give the finite matrix entries explicitly. The orbit of $p$ has cardinality
\[
 \frac{2^{\ell(p)}(N)_{\ell(p)}}{\prod_{j\ge1}m_j(p)!},
 \qquad (N)_\ell=N(N-1)\cdots(N-\ell+1).
\]
For adjacent orbits $p,r$, if $n_{pr}$ is the number of neighbors in the orbit
$r$ of a point in the orbit $p$, the symmetrized adjacency entry is
$\sqrt{n_{pr}n_{rp}}$. Appending a part $1$ gives
\[
 n_{pr}=2(N-\ell(p)),\qquad n_{rp}=m_1(p)+1.
\]
Replacing a part $j\ge1$ by $j+1$ gives
\[
 n_{pr}=m_j(p),\qquad n_{rp}=m_{j+1}(p)+1.
\]
These cases and their transposes exhaust the edges.

Put $z=w^2$, and abbreviate
\[
 Q_{pr}(z)=
 \sqrt{\frac{s(p)s(r)}{(s(p)-\rho(z))(s(r)-\rho(z))}}.
\]
Choose every square root holomorphically and positively at $z=0$. By
\eqref{eff:disk-bounds}, the factors $s(p)-\rho(z)$ have positive real part
on $|z|\le1/25$, so these choices are well defined. The finite matrix
$\mathsf B(w)$ has zero diagonal and the following entries:
\begin{equation}\label{eff:complex-orbit-entries}
 \mathsf B_{pr}(w)=
 \begin{cases}
 w\sqrt{\dfrac{m_1(p)+1}{2}}\,
     \sqrt{1-\ell(p)w^2}\,Q_{pr}(w^2),
      &r=p\cup\{1\},\\[4pt]
 \dfrac{w^2}{2}\sqrt{m_j(p)(m_{j+1}(p)+1)}\,Q_{pr}(w^2),
      &r=(p\setminus\{j\})\cup\{j+1\},
 \end{cases}
\end{equation}
with symmetric entries on reversed edges. Its source vector is
\begin{equation}\label{eff:complex-source}
 \mathsf f_{(1,1)}(w)=
 2\sqrt{\frac{1-w^2}{2-\rho(w^2)}},\qquad
 \mathsf f_{(2)}(w)=w\sqrt{\frac2{4-\rho(w^2)}},
\end{equation}
and vanishes on all other partitions. For real $w=N^{-1/2}$, these are exactly
the compression of $B$ and the vector $f$ in \eqref{eff:DBf}.

All entries are holomorphic for $|w|\le1/5$ and in a neighborhood of this
closed disk. Let $U$ be diagonal with $U_{pp}=(-1)^{\ell(p)}$. Formula
\eqref{eff:complex-orbit-entries} gives
\[
 \mathsf B(-w)=U\mathsf B(w)U,
 \qquad \mathsf f(-w)=U\mathsf f(w).
\]
Consequently each bilinear moment
$\mathsf f(w)^T\mathsf B(w)^{2j}\mathsf f(w)$ is even in $w$ and is a
holomorphic function of $z=w^2$. Denote the sum of these five functions by
$\mathsf S_4(z)$. In particular,
\begin{equation}\label{eff:S-interpolation}
 \mathsf S_4(1/N)=S_4(N).
\end{equation}
The transpose, rather than conjugate transpose, is essential in this
holomorphic expression. Its absolute value is bounded by the usual Euclidean
norms.

A bound on this finite expression in a complex disk will control its Taylor remainder.

\begin{lemma}\label{eff:complex-bound}
The function
\begin{equation}\label{eff:G-definition}
 G(z)=\delta(z)-\mathsf S_4(z)
\end{equation}
is holomorphic in a neighborhood of $|z|\le1/25$ and satisfies
\begin{equation}\label{eff:G-circle}
 |G(z)|<71\qquad(|z|=1/25).
\end{equation}
\end{lemma}

\begin{proof}
Holomorphy follows from the preceding construction and
\eqref{eff:denominator-bound}. To bound the matrix entries, let $\mathsf B^0$
be the real finite matrix on $\mathcal P$ obtained at $N=25$ and $E=50$,
so that $z=1/25$ and $\rho=1$. On $|z|=1/25$, the denominator contribution
in the ratio of corresponding entries is bounded by $10/9$, because
\[
 \frac{s-1}{|s-\rho(z)|}\le
 \frac{s-1}{s-1-1/10}\le\frac{10}{9}
 \qquad(s\ge2).
\]
For a creation edge, the remaining squared ratio is at most
\[
 \frac{|1-\ell z|}{1-\ell/25}
 \le\frac{25+\ell}{25-\ell}\le\frac{31}{19}
 \qquad(0\le\ell\le6).
\]
The factors $|w|$ and $|w|^2$ agree with their real comparison values. Thus
entrywise,
\[
 |\mathsf B(w)|\le\frac{10}{9}\sqrt{\frac{31}{19}}\,\mathsf B^0.
\]
Entrywise domination by a nonnegative matrix also implies domination of the
Euclidean operator norm: apply the entrywise inequality to $|v|$ for any
vector $v$.

The layer estimate \eqref{eff:layer-bound} applies to $\mathsf B^0$. For
$2\le k\le5$, it gives squared block norms at most $9/50$.
The matrix of block norms has at most two nonzero entries in each row,
so its norm is at most $2\sqrt{9/50}$. Therefore
\begin{equation}\label{eff:complex-B-bound}
 \|\mathsf B(w)\|^2
 \le\frac{3100}{1539}\frac{18}{25}
 =\frac{2232}{1539}<\frac32.
\end{equation}
By \eqref{eff:complex-source} and \eqref{eff:disk-bounds},
\begin{equation}\label{eff:complex-f-bound}
 \|\mathsf f(w)\|^2
 =\left|\frac{4(1-z)}{2-\rho(z)}\right|
  +\left|\frac{2z}{4-\rho(z)}\right|
 \le\frac{208}{45}+\frac4{145}<5.
\end{equation}
Cauchy--Schwarz, applied to the bilinear moments by taking absolute values,
now yields
\[
 |\mathsf S_4(z)|
 \le\sum_{j=0}^4\|\mathsf f(w)\|^2\|\mathsf B(w)\|^{2j}
 <5\sum_{j=0}^4(3/2)^j=\frac{1055}{16}<66.
\]
Finally, $|\delta(z)|<5$ proves \eqref{eff:G-circle}.
\end{proof}

\subsection{Cancellation at zero and the effective error}

\begin{proof}[Proof of Theorem~\ref{eff:main}]
Expanding the rational function \eqref{eff:R-intro} at $N=\infty$ gives
\begin{equation}\label{eff:R-expansion}
 \begin{split}
 R(N)={}&2N-4-\frac{20}{3N}-\frac{253}{9N^2}
 -\frac{28409}{180N^3}\\
 &-\frac{477372781}{453600N^4}+O(N^{-5}).
 \end{split}
\end{equation}
For completeness, the coefficients of $z,z^2,z^3,z^4$ in $\delta(z)-4$
are respectively
\[
 \frac{20}{3},\qquad \frac{20a}{3},\qquad
 \frac{20(a^2+b)}{3},\qquad
 \frac{20(a^3+2ab+bc)}{3},
\]
which give the displayed fractions after substitution.
Proposition~\ref{eff:asym-expansion} gives the same expansion for the actual
constant $C(N)$. Hence $C(N)-R(N)=O(N^{-5})$ as integer $N\to\infty$.
The upper derivative bound in \eqref{eff:F-derivative} then implies
\[
 F_N(R(N))=O(N^{-5}).
\]
By \eqref{eff:finite-split}, \eqref{eff:S-interpolation}, and
\eqref{eff:G-definition},
\begin{equation}\label{eff:F-G-tail}
 F_N(R(N))=G(1/N)-\mathcal R_N.
\end{equation}
Lemma~\ref{eff:tail-bound} therefore proves $G(1/N)=O(N^{-5})$.
Since $G$ is holomorphic at zero, its Taylor coefficients of degrees $0$ through
$4$ must vanish. Indeed, if the first nonzero coefficient had degree
$k<5$, then $N^kG(1/N)$ would tend to a nonzero limit, contradicting this
bound along the positive integer sequence.

It follows that $G(z)/z^5$ has a removable singularity at zero. Apply the
maximum-modulus principle on $|z|\le1/25$, using
Lemma~\ref{eff:complex-bound} on its boundary. We obtain
\begin{equation}\label{eff:G-effective}
 |G(1/N)|\le\frac{71\cdot25^5}{N^5}
 =\frac{693359375}{N^5}\qquad(N\ge50).
\end{equation}
Finally, combine \eqref{eff:scalar-error}, \eqref{eff:F-G-tail},
\eqref{eff:G-effective}, and \eqref{eff:tail-constant}:
\begin{align*}
 |C(N)-R(N)|
 &\le |F_N(R(N))|\\
 &\le |G(1/N)|+\mathcal R_N\\
 &<\left(693359375+\frac{2616003280989}{12500}\right)\frac1{N^5}
 <\frac{10^9}{N^5}.
\end{align*}
This proves \eqref{eff:main-bound}. If $N\ge N_\varepsilon$, then
$10^9/N^5\le\varepsilon$, and the strict inequality above proves the
threshold assertion.
\end{proof}

\begin{remark}\label{eff:scope}
The theorem provides the explicit interval
\[
 R(N)-10^9N^{-5}<C(N)<R(N)+10^9N^{-5}\qquad(N\ge50).
\]
For example, its error bound is $10^{-6}$ at $N=1000$ and $10^{-11}$ at
$N=10000$. The constant $10^9$ has not been optimized. The proof establishes
an absolute error estimate; it does not determine the sign of $C(N)-R(N)$.
Only the fixed finite moment matrix is continued in a complex parameter.
The argument requires no continuation of the lattice itself or of $C(N)$
to noninteger dimensions.
\end{remark}

\section{Effective conditions for stability and blow-up}\label{sec:consequences}

The spectral description now gives explicit criteria for the evolution.
In every integer dimension $N\ge3$, the Green bounds furnish sufficient
conditions for stability and blow-up. For $N\ge50$, the rational function
$R(N)$ determines the regime outside a coupling interval of explicitly
controlled width.

\begin{corollary}\label{cor:effective-dynamics}
For an integer $N\ge50$, set $\varepsilon_N=10^9N^{-5}$. Let $\lambda\ge0$.
\begin{enumerate}
\item If $\lambda\le R(N)-\varepsilon_N$, then the semigroups
$S_{\lambda,h}(t)$ are contractive for every $h>0$ and $t\ge0$, and their constant extensions converge to the continuum Friedrichs evolution as in Theorem~\ref{cont:subcritical}.
\item If $\lambda\ge R(N)+\varepsilon_N$, then the mesh-refinement
blow-up and the logarithmic asymptotic of Theorem~\ref{thm:main-bg} hold
for every nonnegative, nonzero $f\in C_c(\R^N\setminus\{0\})$.
\end{enumerate}
More generally, every positive trial vector $f$ in Proposition~\ref{gr:envelopes}
gives the sufficient conditions $\lambda\le1/s_N(f)$ for contractivity
and $\lambda>1/q_N(f)$ for supercritical blow-up.
\end{corollary}
\begin{proof}
Theorem~\ref{eff:main} gives
$R(N)-\varepsilon_N<C(N)<R(N)+\varepsilon_N$.
Thus the first hypothesis implies $\lambda<C(N)$, and the second
implies $\lambda>C(N)$. Apply Theorems~\ref{thm:main-bg} and~\ref{cont:subcritical}.
For the last assertion, use
$1/s_N(f)\le C(N)\le1/q_N(f)$ from
Proposition~\ref{gr:envelopes}, with the strict inequality in the supercritical
condition as stated.
\end{proof}

The length of the rational decision interval is $2\varepsilon_N$.
The explicit majorant is intended to make approximation effective as
the dimension grows. It does not quantify the much smaller errors that
may occur in particular numerical evaluations. For an individual
dimension, the Green bounds and positive physical profiles provide
additional estimates using the same variational threshold.

The scalar norm identity $C(N)=\|T_N\|^{-1}$ and the rational expression
$R(N)$ have different roles. The first is exact for every $N\ge3$; the
second can be evaluated with finitely many rational operations and has
the uniform error proved above. The first nonzero discrepancy between
them is identified in Appendix~\ref{eff:asym-appendix}.

\appendix
\section{Analytic identification of the large-dimension coefficients}
\label{eff:asym-appendix}

The first three terms of the expansion recover \cite[Theorem~1.1]{LizamaAsymptotics}; the subsequent coefficients supply the additional cancellations and the leading error of the rational approximation.

The cancellations used in the effective error estimate concern the sharp
constant on the infinite lattice.  We identify them by bracketing that
constant between two eigenvalues of fixed finite matrices.  In particular,
no convergence assertion about an infinite formal expansion in the dimension
is needed.

\begin{proposition}
\label{eff:asym-expansion}\label{eff:asymptotic-prefix}
As $N$ tends to infinity through the positive integers,
\begin{align}
 C(N)={}&2N-4-\frac{20}{3N}-\frac{253}{9N^2}
 -\frac{28409}{180N^3}-\frac{477372781}{453600N^4}
 \notag\\[-2pt]
 &\hspace{18mm}-\frac{2981659841639}{381024000N^5}
 +O(N^{-6}).
 \label{eff:asym-expansion-eq}
\end{align}
\end{proposition}

\begin{proof}
\emph{Step 1: an exterior lower bound.}
Write $s=s(x)=|x|^2$ and set $h(x)=s(x)^{-2}$ on $\Om$, with
$h(0)=0$.  For $s>1$ define
\[
 t_i=\frac{2x_i}{s+1},\qquad
 T=\frac{4s}{(s+1)^2}<1.
\]
Then $t_i^2\le T$ and $\sum_i t_i^2=T$.  The elementary identity
\[
 (1+t)^{-2}+(1-t)^{-2}-2
 =\frac{2t^2(3-t^2)}{(1-t^2)^2}
 \le \frac{6t^2}{(1-T)^2}
 \qquad(t^2\le T)
\]
gives
\begin{align*}
 s\frac{L_Nh(x)}{h(x)}
 &=s\left[2N-\frac{s^2}{(s+1)^2}
       \sum_{i=1}^N\bigl((1+t_i)^{-2}+(1-t_i)^{-2}\bigr)\right]\\
 &\ge \frac{2Ns(2s+1)}{(s+1)^2}
       -24\left(\frac{s}{s-1}\right)^4.
\end{align*}
The first coefficient is increasing in $s>1$, and the subtracted
quantity is decreasing.  Therefore, for $s\ge9$ and $N\ge28$,
\begin{equation}
 s\frac{L_Nh(x)}{h(x)}
 \ge \frac{171}{50}N-\frac{19683}{512}\ge2N.
 \label{eff:asym-exterior}
\end{equation}

\emph{Step 2: a finite upper and lower bracket.}
Let
\[
 \mathcal F_N=\{x\in\mathbb Z^N:1\le |x|_1\le8\}.
\]
On the finite-dimensional space with norm
$\|u\|_{s^{-1}}^2=\sum_{x\in\mathcal F_N}|u(x)|^2/s(x)$,
let $\mathsf A_N$ be the Dirichlet restriction of $sL_N$ to
$\mathcal F_N$.  Define the lower comparison matrix by
\begin{equation}
 \mathsf B_N=\mathsf A_N-\operatorname{diag}
 \left(s(x)\sum_{\substack{y\sim x\,,\ y\notin\mathcal F_N\\ y\ne0}}
                      \frac{h(y)}{h(x)}\right).
 \label{eff:asym-lower-matrix}
\end{equation}
These are selfadjoint matrices for the displayed weighted inner product.

To justify the lower comparison on the whole lattice, apply the
ground-state identity of Lemma~\ref{lem:gs} in the form
\begin{equation}
 \E[u]=
 \sum_{\substack{\{x,y\}\subset\Om\\x\sim y}}
 h(x)h(y)\left|\frac{u(x)}{h(x)}-\frac{u(y)}{h(y)}\right|^2
 +\sum_{x\in\Om}\frac{L_Nh(x)}{h(x)}|u(x)|^2.
 \label{eff:asym-gs}
\end{equation}
The first sum counts each unoriented edge once.  Retain only its terms
whose two endpoints belong to $\mathcal F_N$.  Their sum, together with
the potential terms on $\mathcal F_N$, is exactly
$\langle\mathsf B_Nu,u\rangle_{s^{-1}}$: expansion of the squares gives
the diagonal
$2N-\sum_{y\notin\mathcal F_N,\,y\ne0}h(y)/h(x)$
in the unweighted energy, and the usual negative adjacency entries.
All discarded edge terms are nonnegative.  For $x\notin\mathcal F_N$
one has $s(x)\ge |x|_1\ge9$, so \eqref{eff:asym-exterior} gives
\[
 \E[u]\ge
 \langle\mathsf B_Nu|_{\mathcal F_N},u|_{\mathcal F_N}\rangle_{s^{-1}}
 +2N\sum_{x\notin\mathcal F_N}\frac{|u(x)|^2}{s(x)}.
\]
Consequently
\[
 C(N)\ge \min\{\lambda_{\min}(\mathsf B_N),2N\}.
\]
Testing $\mathsf B_N$ on a vector supported on $|x|^2=1$ shows
$\lambda_{\min}(\mathsf B_N)\le2N$.  The variational definition of
$C(N)$, tested on functions supported in $\mathcal F_N$, supplies the
upper bound.  Thus
\begin{equation}
 \lambda_{\min}(\mathsf B_N)\le C(N)
 \le\lambda_{\min}(\mathsf A_N),\qquad N\ge28.
 \label{eff:asym-sandwich}
\end{equation}
This lower bound applies to every finitely supported function; it does
not presuppose symmetry of a minimizer on the infinite lattice.

\emph{Step 3: matrices of fixed size and a simple analytic branch.}
Both finite matrices have nonpositive off-diagonal entries on the
connected graph $\mathcal F_N$.  After the weighted unitary conjugation,
a sufficiently large scalar multiple of the identity minus either
matrix is an irreducible nonnegative matrix.  The Perron--Frobenius
theorem therefore makes its lowest eigenvalue simple with a strictly
positive eigenvector.  Uniqueness and normalization imply invariance
of this eigenvector under every signed coordinate permutation.

For $N\ge8$, the orbits in $\mathcal F_N$ are indexed by the fixed set
$\mathcal P_8$ of the $66$ partitions of the integers $1,\ldots,8$.
For a partition $p$ put
\[
 m(p)=\sum_{a\in p}a,\qquad
 \ell(p)=\#p,\qquad s(p)=\sum_{a\in p}a^2.
\]
An invariant function is represented by its actual common value at the
vertices of each orbit, without normalization by orbit size.  In these
coordinates, and with $z=N^{-1}$, the two scaled matrices are exactly
\begin{equation}
 \frac{\mathsf A_N}{2N}=A_0+zA_1,
 \qquad \frac{\mathsf B_N}{2N}=B_0+zB_1.
 \label{eff:asym-affine}
\end{equation}
Here and below these equalities refer to the invariant orbit matrices.
The following rules specify every entry, including all multiplicities.
Initialize the diagonal of $A_0$ at $s(p)$ and all other entries at
zero.  If $p^+$ denotes the partition formed by appending a part $1$,
then, whenever $p^+\in\mathcal P_8$, add
\[
 (A_0)_{p,p^+}=-s(p),\qquad
 (A_1)_{p,p^+}=s(p)\ell(p).
\]
For each distinct part $a$ of multiplicity $r_a(p)$, and each
$\sigma\in\{-1,1\}$, replace one occurrence of $a$ by $a+\sigma$,
deleting it when $a+\sigma=0$, and sort the remaining parts.  If the
resulting nonempty partition $p_{a,\sigma}$ belongs to
$\mathcal P_8$, add
\begin{equation}
 (A_1)_{p,p_{a,\sigma}}\mathrel{+}=-\frac{s(p)r_a(p)}2.
 \label{eff:asym-orbit-entry}
\end{equation}
The empty partition represents the prescribed value zero at the
origin and contributes no entry.  These rules follow because there
are $2(N-\ell(p))$ neighbors obtained from zero coordinates and
$r_a(p)$ neighbors in each of the two directions for a nonzero part
$a$.  Contributions to the same entry are added.

The matrices $B_0,B_1$ have the same off-diagonal entries.  Only rows
with $m(p)=8$ have different diagonals.  In such a row, writing $s=s(p)$,
the changes to the diagonal are
\begin{align}
 (B_0-A_0)_{p,p}&=-\frac{s^3}{(s+1)^2},\notag\\
 (B_1-A_1)_{p,p}&=
 \frac{s^3\ell(p)}{(s+1)^2}
 -\frac{s^3}{2}\sum_{a\in\operatorname{supp}(p)}
                  \frac{r_a(p)}{(s+2a+1)^2}.
 \label{eff:asym-boundary-entries}
\end{align}
Indeed, the first two terms arise from appending a new coordinate,
and the sum arises from increasing an existing coordinate; these
are exactly the edges leaving the core.

Order the partitions by nondecreasing $m(p)$, with $(1)$ first.
Then $A_0$ is upper triangular, with diagonal $s(p)$; its only
nonzero off-diagonal entries append a part $1$.
The matrix $B_0$ is upper triangular as well.  On the boundary
$m(p)=8$, its diagonal is
\[
 d(s)=s-\frac{s^3}{(s+1)^2}
     =\frac{s(2s+1)}{(s+1)^2}\ge\frac{136}{81}>1.
\]
The inequality follows from $s\ge8$ and
$d'(s)=(3s+1)/(s+1)^3>0$.
For either matrix the first diagonal entry is $1$, and every other
nonboundary diagonal entry is at least $2$.  Therefore $1$ is a
simple isolated eigenvalue of both $A_0$ and $B_0$, with right
eigenvector equal to the coordinate vector at $(1)$.

For completeness, analyticity here is a finite-dimensional fact:
apply the analytic implicit function theorem to the characteristic
polynomial at its simple zero $(z,e)=(0,1)$.  A corresponding right
eigenvector, with its $(1)$ coordinate fixed at $1$, is analytic by
solving the remaining equations near that simple eigenvalue.
Denote the resulting eigenvalue branches by $e_A(z)$ and $e_B(z)$.
All other eigenvalues have limits strictly larger than $1$; continuity
of the roots of the characteristic polynomials shows that, for
sufficiently large integer $N$, these two branches are the respective
lowest eigenvalues of the scaled matrices.  Hence
\begin{equation}
 2N e_B(1/N)\le C(N)\le2N e_A(1/N)
 \label{eff:asym-branch-sandwich}
\end{equation}
for all sufficiently large integers $N$.

\emph{Step 4: locality of the coefficient recurrence.}
For either pencil $M_0+zM_1$ in \eqref{eff:asym-affine}, expand
\[
 e(z)=\sum_{k\ge0}e_kz^k,\qquad
 u(z)=\sum_{k\ge0}u_kz^k,
\]
where $e_0=1$, $u_0$ is the coordinate vector at $(1)$, and
$u_k((1))=0$ for $k\ge1$.  Comparing coefficients gives
\begin{equation}
 (M_0-I)u_k=-M_1u_{k-1}
              +\sum_{j=1}^k e_j u_{k-j}.
 \label{eff:asym-recurrence}
\end{equation}
We claim inductively that
\begin{equation}
 u_k(p)=0\quad\text{if }m(p)>k+1.
 \label{eff:asym-support}
\end{equation}
It holds for $k=0$.  The matrix $M_1$ connects only neighboring
levels or acts diagonally, so the right side of
\eqref{eff:asym-recurrence} vanishes at every level exceeding $k+1$,
assuming the claim at earlier orders.  The term $e_ku_0$ is supported
only at the first level.  Solving the upper triangular equations
from the largest level down proves the claim, since each diagonal
entry of $M_0-I$ on those levels is nonzero.

For $k\le6$, \eqref{eff:asym-support} puts $u_k$ on levels at most
$7$.  Both differences $B_0-A_0$ and $B_1-A_1$ are diagonal and
supported on level $8$.  They therefore annihilate $u_k$ and
$u_{k-1}$ in the coefficient equations through order $6$.
The two branches have the same eigenvalue and eigenvector
coefficients through that order.

\emph{Step 5: exact coefficients and the infinite-lattice conclusion.}
We give an explicit recursion for the rational arithmetic involved.
Index the partitions as above by $i=0,\ldots,65$, with $i=0$
corresponding to $(1)$.  Having computed orders below $k$, set
\[
 r_k=-M_1u_{k-1}+\sum_{j=1}^{k-1}e_j u_{k-j},\qquad u_k(0)=0.
\]
For $i=65,64,\ldots,1$ successively, compute
\begin{equation}
 u_k(i)=
 \frac{r_k(i)-\sum_{j>i}(M_0)_{ij}u_k(j)}{(M_0)_{ii}-1},
 \qquad
 e_k=(M_0u_k+M_1u_{k-1})(0).
 \label{eff:asym-exact-recursion}
\end{equation}
The first-row formula can also be written
\begin{equation}
 e_k=-u_k((1,1))+u_{k-1}((1,1))
                        -\tfrac12 u_{k-1}((2)).
 \label{eff:asym-first-row}
\end{equation}
Thus the matrix-entry rules
\eqref{eff:asym-orbit-entry}--\eqref{eff:asym-boundary-entries},
followed by \eqref{eff:asym-exact-recursion}, determine all
coefficients by finite rational operations.
For example, at order one they give
$u_1((1,1))=2$ and $u_1((2))=2/3$, with all higher levels zero,
and consequently $e_1=-2$.  At order two the recurrence gives
$u_2((1,1))=5$ and $u_2((2))=16/9$, so
$e_2=-5+2-1/3=-10/3$.
The successive eigenvalue coefficients are
\[
\begin{array}{c|r}
 k&e_k\\ \hline
 0&1\\
 1&-2\\
 2&-10/3\\
 3&-253/18\\
 4&-28409/360\\
 5&-477372781/907200\\
 6&-2981659841639/762048000
\end{array}
\]
For additional intermediate values useful in checking the recurrence,
\[
\begin{array}{c|r|r}
 k&u_k((1,1))&u_k((2))\\ \hline
 3&109/6&305/36\\
 4&4178/45&4534883/90720\\
 5&89821141/151200&25768127977/76204800\\
 6&61213291711/14112000&1778862208771807/704132352000
\end{array}
\]
Analyticity gives
\[
 e_A(z)=\sum_{k=0}^6e_kz^k+O(z^7),\qquad
 e_B(z)=\sum_{k=0}^6e_kz^k+O(z^7).
\]
Multiplying by $2N$ and applying
\eqref{eff:asym-branch-sandwich} proves
\eqref{eff:asym-expansion-eq} for the actual sharp constant.
\end{proof}

Comparing the next coefficient with the rational expansion determines the leading error of $R(N)$.

\begin{corollary}
\label{eff:asym-leading-error}
Let $R(N)$ be defined by \eqref{eff:R-intro}.  Then
\begin{equation}
 C(N)-R(N)=-\frac{\kappa}{N^5}+O(N^{-6}),\qquad
 \kappa=\frac{250837656967909}{898285248000}>0.
 \label{eff:asym-kappa}
\end{equation}
In particular, $R(N)>C(N)$ for all sufficiently large integers $N$,
and $\lim_{N\to\infty}N^5(R(N)-C(N))=\kappa$.
\end{corollary}
\begin{proof}
With $z=N^{-1}$, expand the rational defect as
\[
 2N-4-R(N)=\frac{20}{3}
 \frac{z(1-cz)}{1-(a+c)z+(ac-b)z^2}
 =\sum_{j\ge1}d_jz^j.
\]
The coefficients are determined by
\[
 d_0=0,\qquad d_1=\frac{20}{3},\qquad
 d_2=(a+c)d_1-\frac{20c}{3},\qquad
 d_j=(a+c)d_{j-1}-(ac-b)d_{j-2}\quad(j\ge3).
\]
They give
\[
 d_1=\frac{20}{3},\quad d_2=\frac{253}{9},\quad
 d_3=\frac{28409}{180},\quad
 d_4=\frac{477372781}{453600},\quad
 d_5=\frac{61007319607185121}{8084567232000}.
\]
Subtracting this convergent rational expansion from
\eqref{eff:asym-expansion-eq} gives \eqref{eff:asym-kappa}, since
\[
 \frac{2981659841639}{381024000}
 -\frac{61007319607185121}{8084567232000}
 =\frac{250837656967909}{898285248000}.
\]
The last number is positive, so the remainder is eventually smaller
in absolute value than $\kappa/(2N^5)$, which proves the eventual sign.
\end{proof}

\begin{remark}
The last corollary identifies an eventual sign and an exact limiting
error coefficient.  Its asymptotic remainder alone supplies no
specified dimension from which that sign holds.  The explicit uniform
absolute-error estimate proved in the main text provides a separate,
effective conclusion.
\end{remark}

\clearpage
\section*{Declarations}

\paragraph{Funding.}
This work was supported by ANID/FONDECYT Regular, grant No.~1260057.

\paragraph{Competing interests.}
The author declares no competing interests relevant to this work.

\paragraph{AI use declaration.}
ChatGPT (OpenAI) was used interactively to assist in exploring proof
strategies, developing and checking symbolic calculations and verification
scripts, checking bibliographic information, and drafting and revising
the manuscript. The finite-dimensional coefficient calculations were
cross-checked with exact-arithmetic scripts; the mathematical arguments
and estimates supporting the results are presented in full in the paper.
The author takes full responsibility for the mathematical content,
the references and the final manuscript.

\end{document}